\documentclass[reqno]{amsart}
\usepackage{hyperref}
\usepackage{mathrsfs}

\AtBeginDocument{}

\begin{document}
\title[\hfil Regularity of pullback attractors and equilibrium]
{ Regularity of pullback attractors and equilibrium for  non-autonomous stochastic
 FitzHugh-Nagumo system on unbounded domains}

\author{Wenqiang Zhao, Anhui Gu}  

\address{Wenqiang Zhao \newline
School Of Mathematics and Statistics, Chongqing Technology and Business
University,  Chongqing 400067, China}
\email{gshzhao@sina.com}

\address{Anhu Gu \newline
 School of Mathematics  and statistics, Southwest University,
Chongqing 400715, China}
\email{gahui@swu.edu.cn}


\begin{abstract}
 A theory on bi-spatial random attractors developed recently  by Li \emph{et al.} is extended to study  stochastic Fitzhugh-Nagumo system driven by a non-autonomous
term as well as a general multiplicative noise. By using the so-called  notions of  uniform absorption and uniformly pullback asymptotic compactness, it is showed that
every generated random cocycle has a pullback attractor in $L^l(\mathbb{R}^N)\times L^2(\mathbb{R}^N)$ with $l\in(2,p]$, and the family of obtained attractors is upper
 semi-continuous at any intensity of noise. Moreover,  if some additional conditions are added,  then the system possesses a unique equilibrium and  is attracted by
 a single point.
\end{abstract}

\maketitle
\numberwithin{equation}{section}
\newtheorem{theorem}{Theorem}[section]
\newtheorem{lemma}[theorem]{Lemma}
\newtheorem{remark}[theorem]{Remark}
\newtheorem{definition}[theorem]{Definition}
\allowdisplaybreaks

\section{Introduction}

In this paper, we consider the random dynamics of solutions  of
the following non-autonomous FitzHugh-Nagumo system defined on $\mathbb{R}^N$ perturbed by a $\varepsilon$-multiplicative noise:
\begin{align}\label{FN1}
d\tilde{u}+(\lambda\tilde{ u}-\Delta\tilde{ u}+\alpha \tilde{v})dt=f(x,\tilde{u})dt+g(t,x)dt+\varepsilon \tilde{u}\circ d\omega(t),t> \tau,
 \end{align}
\begin{align}\label{FN2}
d\tilde{v}+(\sigma \tilde{v}-\beta\tilde{u})dt=h(t,x)dt+\varepsilon\tilde{v}\circ d\omega(t),t>\tau,
 \end{align}
with initial value numbers
\begin{align}\label{FN3}
\tilde{u}(x,\tau)=\tilde{u}_\tau(x), \ \ \  \ \ \ \tilde{v}(x,\tau)=\tilde{v}_\tau(x),
\end{align}
where the initial condition $(\tilde{u}_0,\tilde{v}_0)\in L^2(\mathbb{R}^N)\times L^2(\mathbb{R}^N)$, the coefficients $\lambda,\alpha,\beta,\sigma$ are  positive constants,
the non-autonomous terms $g,h\in L^2_{loc}(\mathbb{R}, L^2(\mathbb{R}^N))$, the nonlinearity $f$ is a smooth function satisfying some polynomial growth, $\varepsilon$ is the intensity of noise with $\varepsilon\in [-a,a]\setminus \{0\},a>0$, $\omega(t)$ is a Wiener process defined on a probability space  $(\Omega,\mathcal{F},P)$, where
 $\Omega=\{\omega\in
C(\mathbb{R},\mathbb{R}); \omega(0) =0\}$, and
 $\mathcal {F}$ be the
 Borel $\sigma$-algebra induced by the compact-open topology of
$\Omega$ and ${P}$ be  the corresponding
 Wiener measure on $(\Omega,\mathcal{F})$.

The deterministic FitzHugh-Nagumo system is an important mathematical model to describe the signal transmission across axons in neurobiology, see \cite{Nagumo,FitzHugh,Jones,Bell} and references therein.   It is well studied in the
literature, see e.g. \cite{Sys1,Sys2,Sys3,Nagumo1}. In the random case,
when $g$ and $h$ do not depend on the time, Wang \cite{Wang2} proved the existence and uniqueness of random attractors in $L^2(\mathbb{R}^N)\times L^2(\mathbb{R}^N)$.
For the general non-autonomous forcings $g$ and $h$,  under
additive noises,  Adili and Wang \cite{Adili2} obtained  the pullback attractors in $L^2(\mathbb{R}^N)\times L^2(\mathbb{R}^N)$, and
Bao \cite{Bao} developed  this result and obtained the regularity in  $H^1(\mathbb{R}^N)\times L^2(\mathbb{R}^N)$.
For our problem (\ref{FN1})-(\ref{FN3}) ,  \emph{i.e.}, under multiplicative noise, Adili and Wang \cite{Adili1} proved  the existence and upper semi-continuity of attractors
in $L^2(\mathbb{R}^N)\times L^2(\mathbb{R}^N)$ recently. For the stochastic lattice FitzHugh-Nagumo system, the existences of random attractors
 are widely studied  in \cite{Huang,Gu,Gu1}.  However, to our knowledge, there are no literature to investigate the asymptotic high-order integrability of solutions to
the FitzHugh-Nagumo system, even for the deterministic case.

In this paper,  we strengthen these results offered by  \cite{Adili1} and  devote to obtain the asymptotic high-order integrability of solutions of problem (\ref{FN1})-(\ref{FN3}).
To this end, a theory on bi-spatial random attractors developed recently  by Li \emph{et al.} \cite{Liyangrong0, Liyangrong1} is extended  to stochastic partial differential equations (SPDE)  with both non-autonomous terms and
random noises, see Theorem 2.9 and 2.10.  It is showed that the \emph{uniform absorption} and \emph{uniformly pullback asymptotic compactness} are
the appropriate notions to depict the existence and upper semi-continuity of attractor in both the initial space and terminate space \cite{Liyangrong0}.
As for the theory on the upper semi-continuity of pullback attractors in an initial space and its applications, we may refer to \cite{Wang3,Wang1,Wang5,Zhou1} and references therein.

Then we apply the obtained theorems to prove that  the problem (\ref{FN1})-(\ref{FN3})  admits  a unique pullback
attractor in $L^p(\mathbb{R}^N)\times L^2(\mathbb{R}^N)$, with the functions $f,g $ and $h$  satisfying almost the same conditions  as in \cite{Adili1}.
Furthermore, we derive the upper semi-continuity of pullback attractors of system in $L^p(\mathbb{R}^N)\times L^2(\mathbb{R}^N)$ as the intensity
 $\varepsilon$ approaches any $\varepsilon_0\in\mathbb{R}$. These are achieved by checking the uniform absorption and uniformly pullback asymptotic compactness properties of random cocycles.
 The uniformly pullback asymptotic compactness in $L^p(\mathbb{R}^N)\times L^2(\mathbb{R}^N)$ is completely proved by estimate of the $L^2$ and $L^p$-uniform boundedness and the $L^p$-truncation of solutions, see Lemma 4.4 and 4.5. it seems that the estimate of  $L^2$-truncation is unnecessary, see
 \cite{Liyangrong0,Liyangrong1,Liyangrong2,Zhao3,Zhao4,Lijia,Yin,Zhong}.
 It is worth mention that an additional assumption on the non-autonomous terms (see \cite{Adili1}) is not used in our proof, see section 3.

The third goal is to study \emph{stochastic fixed point} or \emph{random equilibrium} of random dynamical system, see \cite{Chues}. In this paper, we introduce the notion of equilibrium for SPDE with both non-autonomous terms and white noises. It is showed that if the parameters satisfy some additional conditions, then the system admits a unique
 equilibrium and  is attracted by a single point.

This paper is organized as follows. In the next section, we introduce some concepts required for our further discussions and extend the results developed by
\cite{Liyangrong0,Liyangrong1} to the general SPDE with non-autonomous forcing. In section 3, we give the assumptions on $g, h$ and $f$, and define a family of continuous random cocycles for problem (\ref{FN1})-(\ref{FN3}). In section 4,  we prove the existence and upper semi-continuity of pullback attractors in $L^p(\mathbb{R}^N)\times L^2(\mathbb{R}^N)$.
The final section is concerned with the existence of equilibrium of the random cocycle.

\section{Preliminaries and abstract results}

In this section, we give the sufficient conditions for the existence and upper semi-continuity of pullback attractors in \emph{the terminate space} for random dynamical systems over two parametric spaces, which are applicable to SPDE with both non-autonomous deterministic and random terms, and the structure of the pullback attractor is presented. This is an extension of the
 corresponding results just  established by Li \emph{et al.} \cite{Liyangrong0}, which is suitable for random systems  with only one stochastic terms. The reader is also referred to
\cite{Wang3,Wang4,Zhou2} for the theory of pullback attractors and its applications in \emph{the initial space} over two parametric spaces, and to \cite{Rand1,Rand2,Rand3,Rand4} for  one parametric random attractors. The reader may also refer to \cite{James,Robin} for pullback attractor of deterministic dynamical systems.

\subsection{Preliminaries}

Let both $(X,\|.\|_X)$ and $(Y,\|.\|_Y)$ be separable Banach spaces, where $X$ is called \emph{an initial space} which contains all initial data, and  $Y$ is called \emph{a terminate space} which contains all values of solutions for a SPDE \cite{Liyangrong0}. Both $X$ and $Y$ may not be embedded in any direction, but we assume that they are \emph{limit-uniqueness} in the following sense:\\

\emph{\textbf{(H1)}\ \  If $\{x_n\}_n\subset X\cap Y$ such that $x_n\rightarrow x$ in
$X$ and $x_n\rightarrow y$ in $Y$, respectively, then we have $x=y$.}\\

Let $Q$ be a nonempty set and $(\Omega,\mathcal{F},{P})$ be a probability space. We assume that there are two groups $\{\sigma_t\}_{t\in\mathbb{R}}$ and
$\{\vartheta_t\}_{t\in\mathbb{R}}$ over $Q$ and $\Omega$, respectively. Specifically,  the mapping $\sigma: \mathbb{R}\times Q\mapsto Q$ satisfies that
$\sigma_0$ is the identity on $Q$, and $\sigma_{s+t}=\sigma_s\circ \sigma_t$ for all $s,t\in\mathbb{R}$. Similarly,  $\vartheta: \mathbb{R}\times \Omega\mapsto \Omega$
is a  $(\mathcal{B}(\mathbb{R})\times\mathcal{F}, \mathcal{F})$-measurable  mapping such that $\vartheta_0$ is the identity on $\Omega$, $\vartheta_{s+t}=\vartheta_s\circ \vartheta_t$ for all $s,t\in\mathbb{R}$ and $\vartheta_tP=P$ for all $t\in\mathbb{R}$. In particular, we call both $(Q,\{\sigma_t\}_{t\in\mathbb{R}})$ and
$(\Omega,\mathcal{F},{P},\{\vartheta_t\}_{t\in\mathbb{R}})$  parametric dynamical systems. Let $\mathbb{R}^+=\{x\in\mathbb{R}; x\geq0\}$ and $2^X$ be the collection of
all subsets of $X$.\\

\textbf{Definition 2.1.}\ \ \emph{A measurable mapping $\varphi: \mathbb{R}^+\times Q\times \Omega\times X \rightarrow X, (t,q,\omega,x)\mapsto \varphi(t,q,\omega,x)$ is called
to be a random cocycle on $X$ over $(Q,\{\sigma_t\}_{t\in\mathbb{R}})$ and
$(\Omega,\mathcal{F},{P},\{\vartheta_t\}_{t\in\mathbb{R}})$  if for every $q\in Q,\omega\in \Omega$ and $s,t\in\mathbb{R}^+ $ the following statements are satisfied:
\begin{align}
&(i)\ \varphi(0, q,\omega,.)\ \mbox{is the identity on X}, \ \ \  (ii)\ \varphi(t+s,q, \omega,.)=\varphi(t,\sigma_sq, \vartheta_s\omega, .)\circ \varphi(s,q, \omega, .).\notag
 \end{align}
An random cocycle is said to be continuous in $X$ iff each operator $\varphi(t,q,\omega,.)$ is continuous in $X$ for each  $q\in Q, \omega\in \Omega$ and $t\in\mathbb{R}^+ $.\\ }

In particular, it is pointed out that  in this paper, we need further to assume that the random cocycle $\varphi$ acting on $X$  takes its values  in  the terminate space $Y$ for for all $t>0$ (except that $t=0$), \emph{i.e.,}\\

\emph{\textbf{(H2)}\ \   For every $t>0,q\in Q,$ and $\omega\in\Omega$, $\varphi(t,q,\omega,.): X\rightarrow  Y$.}\\

In the sequel, we use  $\mathcal{D}$ to denote  a collection of some families of nonempty subsets  of $X$ which is parameterized by $(q,\omega)\in(Q\times\Omega)$:
$$
\mathcal{D}=\{D=\{\emptyset\neq D(q,\omega)\in 2^X; q\in Q, \omega\in\Omega \};f_D\ \mbox{satisfies some conditions}\}.
$$
We further assume that $\mathcal{D}$ is inclusion closed, that is, for each $D\in\mathcal{D}$,
$$
\{\tilde{D}(q,\omega): \tilde{D}(q,\omega)\ \mbox{is a nonempty subset of}\ D(q,\omega), \forall\ q\in Q, \omega\in \Omega \}\in\mathcal{D}.
$$
Given $D_1, D_2\in\mathcal{D}$, we say that $D_1=D_2$ if and only if $D_1(q,\omega)=D_2(q,\omega)$ for each $q\in Q$ and $\omega\in \Omega$.

Throughout this paper, all assertions about $\omega$ are assumed to hold on a $\vartheta_t$-invariant set of full measure
(unless some exceptional cases  needed). \\

\textbf{Definition 2.2.}\emph{ A set-valued mapping $K: Q\times \Omega\rightarrow 2^X$ is called measurable in $X$ with respect to $\mathcal{F}$ in $\Omega$ if
 the mapping $\omega\in\Omega\mapsto \mbox{dist}_X(x,K(q,\omega))$ is ($\mathcal{F},\mathcal{B}(\mathbb{R}))$-measurable for every fixed $x\in X$
and $q\in Q$, where $\mbox{dist}_X$ is the Haustorff semi-metric in $X$, i.e., for the two nonempty subsets $A,B\in 2^X$,
$$
dis_X(A,B)=\sup_{a\in A}\inf_{b\in B}\|a-b\|_X.
$$}

\textbf{Definition 2.3.} \emph{ Suppose $\varphi$ is a random cocycle on $X$  over $(Q,\{\sigma_t\}_{t\in\mathbb{R}})$ and
$(\Omega,\mathcal{F},{P},\{\vartheta_t\}_{t\in\mathbb{R}})$ and takes its value in $Y$. A set valued mappping $\mathcal{A}: Q\times\Omega\mapsto 2^{X\cap Y}$ is called  a
(X,Y)-pullback attractor for  $\varphi$  if}

\emph{(i)  $\mathcal {A}$ is measurable in $X$  (w.r.t the $P$-completion  of  $\mathcal{F}$ in $\Omega$),  and
$\mathcal {A}(q,\omega)$ is compact in $Y$ for all $q\in Q, \omega\in\Omega$,}

 \emph{(ii) $\mathcal {A}$ is invariant, that is, for every $q\in Q, \omega\in\Omega$,
 $$
 \varphi(t,q,\omega,\mathcal{A}(q,\omega))=\mathcal{A}(\sigma_tq,\vartheta_t\omega), \forall\ t\geq0,
 $$}

\emph{(iii) $\mathcal {A}$  attracts every element $D\in \mathcal{D}$ in $Y$, that is, for every $q\in Q, \omega\in\Omega$,
$$\mathop {\lim }\limits_{t\rightarrow +\infty}\mbox{dist}_Y(\varphi(t,\sigma_{-t}q,\vartheta_{-t}\omega, D(\sigma_{-t}q, \vartheta_{-t}\omega)),\mathcal{A}(q,\omega))= 0,$$
}
where $dist_Y$ is the Haustorff semi-distance in $Y$ and the set $\varphi(t,q,\omega,D(q,\omega))=\{\varphi(t,q,\omega,x);x\in D(q,\omega)\}$.
\\

If $X=Y$, the above concept reduces to the well known notion of a $\mathcal{D}$-pullback attractors, which is first introduced in \cite{Wang4}. We also remark that the
measurability of $\mathcal{A}$ is assumed in the initial space $X$.\\

\textbf{Definition 2.4.}(see\cite{Jones}.)\ \  \emph{ Let both $Z$ and $I$ be two metric spaces. A family $\{\mathcal{A}_{\alpha}\}_{\alpha\in I}$ of sets in $Z$ is said to be upper semi-continuous at $\alpha_0$
if
\begin{align}\label{}
\lim_{\alpha\rightarrow\alpha_0}\mbox{dist}_Z(\mathcal{A}_\alpha,\mathcal{A}_{\alpha_0})=0.\notag
\end{align}
A family $\mathcal{A}_{\alpha}$ of set-mappings over $Q$ and $\Omega$ is called to be upper semi-continuous if $\mathcal{A}_{\alpha}(q,\omega)$ is upper semi-continuous for each
$q\in Q$ and $\omega\in \Omega$.}\\

In the sequel, we need to consider a family of random cocycles $\{\varphi_\alpha\}_{\alpha\in I}$ with $I=[-a,a]\setminus\{0\}$, where $a>0$ and  $\varphi_0$ is a deterministic cocycle over the parametric space $(Q,\{\sigma_t\}_{t\in\mathbb{R}})$.\\

\textbf{Definition 2.5.}(see\cite{Liyangrong0}.) \ \ \emph{ A family of random cocycles $\{\varphi_\alpha\}_{\alpha\in I}$ is said to be convergent at point $\alpha=\alpha_0$ in $X$ if for each $q\in Q, \omega\in\Omega,$ and $x,x_0\in X$,
$$
\varphi_\alpha(t, q, \omega,x)\rightarrow \varphi_{\alpha_0}(t, q, \omega,x_0)\ \mbox{in}\ X,
$$
 whenever $\alpha\rightarrow\alpha_0$ and $x\rightarrow x_0$.  A family of random cocycles $\{\varphi_\alpha\}_{\alpha\in I}$ is said to be convergent in $X$ if it is convergent at
any point $\alpha$. We say a family of random cocycles  $\alpha_\varepsilon(\alpha\in(0,a])$ converges to a deterministic cocycle $\varphi_0$ in $X$ if for each $q\in Q,$ and $x,x_0\in X$,
$$
\varphi_\alpha(t, q, \omega,x)\rightarrow \varphi_{0}(t, q ,x_0)\ \mbox{in}\ X,
$$
 whenever $\alpha\rightarrow\alpha_0$ and $x\rightarrow x_0$.} \\

\textbf{Definition 2.6.}(see\cite{Liyangrong0}.) \ \ \emph{ A family of random cocycles $\{\varphi_\alpha\}_{\alpha\in I}$ is said to be uniformly absorbing in $X$ if
each $\varphi_\alpha$  has a closed and measurable pullback absorbing set $K_\alpha$ in $X$ such that the closure $\overline{K}=\{\overline{\cup_{\alpha\in I}K_\alpha(q,\omega)};q\in Q, \omega\in\Omega\} \in\mathcal{D}$ and for each $q\in Q, \omega\in\Omega$,
$$
\limsup_{\alpha\rightarrow0}\|K_\alpha(q,\omega)\|_X\leq c\  \ \mbox{for some deterministic constant}\ c>0.
$$
Here a pullback absorbing set $K_\alpha$ means that for every $D\in \mathcal{D}$, there exists an absorbing time $T=T(D,q,\omega)>0$ such that for each
$q\in Q, \omega\in\Omega$,
$$\varphi_\alpha(t,\sigma_{-t}q, \vartheta_{-t}\omega, D(\sigma_{-t}q,\vartheta_{-t}\omega))\subseteq K(q, \omega)\ \ \ \ \ \ \mbox{for all}\ t\geq
T.$$}

\textbf{Definition 2.7.} \ \emph{A family of random cocycles $\{\varphi_\alpha\}_{\alpha\in I}$ is said to be uniformly pullback asymptotically compact over $I$ in $X$  if for each $q\in Q, \omega\in\Omega$, $D\in \mathcal{D}$, the sequence
\begin{align}\label{uniform}
\{\varphi_{\alpha_n}(t_n, \sigma_{-t_n}q, \vartheta_{-t_n}\omega, x_n)\}\ \mbox{has a convergence subsequence in}\ X,
\end{align}
whenever $\alpha_n\in I, t_n\rightarrow \infty$, and
$x_n\in D(\sigma_{-t_n}q,\vartheta_{-t_n}\omega)$. A family of random cocycles $\{\varphi_\alpha\}_{\alpha\in I}$  is uniformly pullback asymptotically compact in $Y$ if the convergence in (\ref{uniform}) holds under $Y$-norm. A single cocycle $\varphi_{\alpha_0}$ is pullback asymptotically compact in $X$ if (\ref{uniform}) holds for a single point $\alpha=\alpha_0$.}  \\

\textbf{Definition 2.8.} \emph{Let $\mathcal{D}$ be a collection of some families of nonempty subsets of $X$, and $\varphi$ be a random cocycle on $X$ over
$(Q,\{\sigma_t\}_{t\in\mathbb{R}})$ and $(\Omega,\mathcal {F},{P},\{\vartheta_t\}_{t\in\mathbb{R}})$. A mapping $\psi: \mathbb{R}\times Q\times \Omega\rightarrow X$ is called \emph{a complete orbit}
of $\varphi$ if for each $\tau\in\mathbb{R}, t\in\mathbb{R}^+,q\in Q$ and $\omega\in \Omega$, there holds:
$$
\varphi(t,\sigma_{\tau}q,\sigma_{\tau}\omega,\psi(\tau,q,\omega))=\psi(t+\tau,q,\omega).
$$ If  in addition, there exists $D=\{D(q,\omega);q\in Q,\omega\in \Omega\in\mathcal{D}\}$ such that $\psi(\tau,q,\omega)\in D(\sigma_{\tau}q,\sigma_{\tau}\omega)$, then $\psi$
is called a $\mathcal{D}$-complete orbit of $\varphi$.}

\subsection{Abstract results}

By slightly modifying the arguments of Theorem  3.1 in  Li \emph{et al}. \cite{Liyangrong0}, we can extend the corresponding theory to the case  in which the random  systems
with non-autonomous term as well as random noises are considered.\\

 \textbf{ Theorem  2.9.}  \emph{Let $(X,Y)$ be a pair of Banach spaces satisfying hypothesis \emph{(H1)}, and $\varphi$  a continuous random cocycle  in $X$ over $(Q,\{\sigma_t\}_{t\in\mathbb{R}})$ and
$(\Omega,\mathcal{F},{P},\{\vartheta_t\}_{t\in\mathbb{R}})$  such that hypothesis \emph{(H2)} holds.
 Assume further that\\}

\emph{(i) $\varphi$ has  a closed and measurable (w.r.t. the $P$-completion of $\mathcal{F}$) pullback absorbing set
$K=\{K(q,\omega);q\in Q,\omega\in\Omega\}\in \mathcal{D}$ in $X$;}

\emph{(ii) $\varphi$ is pullback asymptotically compact in $X$;}

\emph{(iii) $\varphi$ is pullback asymptotically compact in $Y$.}\\

\emph{Then the random cocycle $\varphi$ admits a unique (X,Y)-pullback attractor $\mathcal{A}\in \mathcal{D}$, which is  structured by
\begin{align}\label{AAAA}
\mathcal{A}(q,\omega)&=\cap_{\tau>0}\overline{\cup_{t\geq \tau}(t,\sigma_{-t}q,\vartheta_{-t}\omega, K(\sigma_{-t}q,\vartheta_{-t}\omega))}^Y\notag\\
&=\{\psi(0,q,\omega); \psi \ \mbox{is a}\ \mathcal{D}\mbox{-complete orbit of the random cocycle}\ \varphi\}.
 \end{align}
Moreover, $\mathcal{A}=\mathcal{A}_X$, where $\mathcal{A}_X$ is the $(X,X)$-pullback attractor.}\\

In the following,  we will consider both the semi-continuity of a family of  the bi-spatial pullback attractors and the existence problem. We give a unified result, where the concepts of uniform
absorption, uniform convergence  and uniformly pullback asymptotic compactness are used.  To this end, we need to consider a family of random cocycles $\{\varphi_\alpha\}_{\alpha\in I}$ with $I=[-a,a]\setminus\{0\}$, where $a>0$ and  $\varphi_0$ is a deterministic cocycle over the parametric space $(Q,\{\sigma_t\}_{t\in\mathbb{R}})$.

Then analogous to
 the proofs of Theorem 3.1  in
Li \emph{et al.} \cite{Liyangrong1}, or some small modifying the argument of   Theorem 4.1 in Li \emph{et al.} \cite{Liyangrong0}, we can extend their results about random attractors
 with single random noises to the random cocycle  with non-autonomous term as well as random noises. \\

\textbf{ Theorem  2.10.}  Let $\mathcal{D}$ be a collection of some families of
 nonempty subsets of $X$ and $(X,Y)$ a pair of Banach spaces satisfying hypothesis \emph{(H1)}. Suppose that $\{\varphi_\alpha\}_{\alpha\in I}$ is a family of continuous random cocycles in $X$ over  $(Q,\{\sigma_t\}_{t\in\mathbb{R}})$ and
$(\Omega,\mathcal{F},{P},\{\vartheta_t\}_{t\in\mathbb{R}})$  such that hypothesis \emph{(H2)} holds,  and  $\varphi_0$ is  a continuous deterministic cocycle in $X$  over $(Q,\{\sigma_t\}_{t\in\mathbb{R}})$  satisfying
$\varphi_0(t,q,.): X\rightarrow Y$ for all $t>0$ and $q\in Q$.
 Assume further that\\

\emph{(i) $\varphi_\alpha$  is  convergent in $X$ at any $\alpha\in [-a,a]$;}

\emph{(ii)  $\varphi_\alpha(\alpha\in I)$  is  uniformly absorbing in  $X$ ;}

\emph{(iii) $\varphi_\alpha(\alpha\in I)$ is uniformly pullback asymptotically compact in $X$;}

\emph{(iv) $\varphi_\alpha(\alpha\in I)$ is uniformly pullback asymptotically compact in $Y$.}\\

\emph{Then each  random cocycle $\varphi_\alpha$ admits a unique (X,Y)-pullback attractor $\mathcal{A}_\alpha\in \mathcal{D}$, such that the family
 $\{\varphi_\alpha\}_{\alpha\in I}$ is upper semi-continuous at any $\alpha\in I$ in both $X$ and $Y$.  If in addition, $\varphi_0 $ has an (X,Y)-attracting set $A_0$,
 then the $\{\varphi_\alpha\}_{\alpha\in I}$ is upper semi-continuous at $\alpha=0$ in both $X$ and $Y$.}

\section{ Non-autonomous FitzHugh-Nagumo system on $\mathbb{R}^N$ with multiplicative noise}

For the non-autonomous FitzHugh-Nagumo system (\ref{FN1})-(\ref{FN3}), the nonlinearity $f(x,s)$
satisfy almost the same assumptions as in \cite{Adili1}, \emph{i.e.},
for $ x\in\mathbb{R}^N$ and $s\in\mathbb{R}$,
\begin{align}\label {a1}
&f(x,s)s\leq -\alpha_1|s|^p+\psi_1(x),
\end{align}
\begin{align}\label {a2}
&|f(x,s)|\leq \alpha_2|s|^{p-1}+\psi_2(x),
\end{align}
\begin{align}\label {a3}
&\frac{\partial f}{\partial s}f(x,s)\leq \alpha_3 ,
\end{align}
\begin{align}\label {a4}
&\Big|\frac{\partial f}{\partial x}f(x,s)\Big|\leq \psi_3(x).
\end{align} where $p\geq2$, $\alpha_i> 0(i=1,2,3)$ are determined constants,  $\psi_1\in L^{1}(\mathbb{R}^N)\cap L^{\frac{p}{2}}(\mathbb{R}^N)$,
and $\psi_2,\psi_3\in L^2(\mathbb{R}^N)$. The non-autonomous terms
 $g\in L^2_{loc}(\mathbb{R}, L^2(\mathbb{R}^N))$ and $h\in L^2_{loc}(\mathbb{R}, L^2(\mathbb{R}^N))$  satisfy that for every $\tau\in\mathbb{R}$ and some $0<\delta_0<\delta=\min\{\lambda,\sigma\}$,
\begin{align}\label {a5}
\int\limits_{-\infty}^\tau e^{\delta_0 s} (\|g(s,.)\|_{L^2(\mathbb{R}^N)}^2+\|h(s,.)\|_{L^2(\mathbb{R}^N)}^2)ds<+\infty,
\end{align} where $\lambda$ and $\delta$  are as in (\ref{FN1})-(\ref{FN3}). The $H^1$-condition on the non-autonomous term $h$ in (\ref{a5}) is required to prove the asymptotic compactness
of solutions in $L^2(\mathbb{R}^N)\times L^2(\mathbb{R}^N)$, see \cite{Adili1}:
$$
\int\limits_{-\infty}^\tau e^{\delta_0 s} \|h(s,.)\|_{H^1(\mathbb{R}^N)}^2ds<+\infty.
$$

 In order to model the random noises in system (\ref{FN1})-(\ref{FN3}),
  we need to define a shift operator $\vartheta$ on $\Omega$ (which is defined in the introduction) by
$\vartheta_t\omega(s)=\omega(s+t)-\omega(t)$ for every $\omega\in\Omega, t,s\in\mathbb{R}.$
 Then $\vartheta_t$ is a measure preserving transformation group on  $(\Omega,\mathcal{F},P)$, that is, $(\Omega,\mathcal{F},P,\{\vartheta_t\}_{t\in\mathbb{R}})$
 is a parametric dynamical system.  By the law of the iterated logarithm (see \cite{Rand1}), there exists a
 $\vartheta_t$-invarant set $\tilde{\Omega}\subset\Omega$ of full measure sucht that for $\omega
 \in\tilde{\Omega}$,
\begin{align} \label{3.6}
\frac{\omega(t)}{t}\rightarrow0, \ \mbox{as}\ |t|\rightarrow+\infty.
\end{align}
Put $Q=\mathbb{R}$. Define a family of shift operator $\{\sigma_t\}_{t\in\mathbb{R}}$ by $\sigma_t(\tau)=t+\tau$ for all $t,\tau\in\mathbb{R}$. Then both $\{\mathbb{R},\{\sigma_t\}_{t\in\mathbb{R}}\}$ and $(\Omega,\mathcal{F},P,\{\vartheta_t\}_{t\in\mathbb{R}})$ are parametric dynamical systems.  We will define a continuous random
cocycle for system (\ref{FN1})-(\ref{FN3}) over $(Q,\{\sigma_t\}_{t\in\mathbb{R}})$ and $(\Omega,\mathcal{F},P,\{\vartheta_t\}_{t\in\mathbb{R}})$.

Given  $\omega\in \Omega$, put $z(t,\omega)=z_\varepsilon(t,\omega)=e^{-\varepsilon \omega (t)}$. Then we have $dz+\varepsilon z\circ d\omega(t)=0.$
Let $(\tilde{u},\tilde{v})$ satisfy  problem (\ref{FN1})-(\ref{FN3}) and  write
\begin{align} \label{trans}
u(t,\tau,\omega,u_0)=z(t,\omega)\tilde{u}(t,\tau,\omega,\tilde{u}_0)\ \  \mbox{and}\ v(t,\tau,\omega,v_0)=z(t,\omega)\tilde{v}(t,\tau,\omega,\tilde{v}_0).
\end{align}
 Then $(u,v)$ solves the follow system
\begin{align} \label{pr1}
\frac{du}{dt}+\lambda u-\Delta u+\alpha v=z(t,\omega)f(x,z^{-1}(t,\omega)u)+z(t,\omega)g(t,x),
\end{align}
\begin{align} \label{pr2}
\frac{dv}{dt}+\sigma v-\beta u=z(t,\omega)h(t,x),
\end{align}
with initial conditions  $u_\tau=z(\tau,\omega)\tilde{u}_\tau$ and $v_\tau=z(\tau,\omega)\tilde{v}_\tau$.

It is known (see \cite{Adili1}) that for every $(u_\tau,v_\tau)\in L^2(\mathbb{R}^N)\times  L^2(\mathbb{R}^N)$
 the problem (1.1)  possesses  a unique solution $(u,v)$ such that
$u\in C([\tau,+\infty),L^2(\mathbb{R}^N))\cap L^2(\tau,T, H^1(\mathbb{R}^N))\cap L^p(\tau,T, L^p(\mathbb{R}^N))$ and $v\in C([0,+\infty), L^2(\mathbb{R}^N))$. In addition,
the solution $(u,v)$ is continuous in $L^2(\mathbb{R}^N)\times L^2(\mathbb{R}^N)$ with respect to the initial value $(u_\tau,v_\tau)$ . Then formally $(\tilde{u},\tilde{v})=(z^{-1}(t,\omega)u,z^{-1}(t,\omega)v)$ is the solution to
problem (\ref{FN1})-(\ref{FN3}) with the initial value $\tilde{u}_\tau=z^{-1}(\tau,\omega){u}_\tau$ and $\tilde{v}_\tau=z^{-1}(\tau,\omega){v}_\tau$.

We are at the position to give the continuous random cocycle $\varphi$ associated with problem (\ref{FN1})-(\ref{FN3}) over $(Q,\{\sigma_t\}_{t\in\mathbb{R}})$ and $(\Omega,\mathcal{F},P,\{\vartheta_t\}_{t\in\mathbb{R}})$. Define
\begin{align}\label{eq0}
& \ \ \ \ \ \ \ \ \ \varphi (t,\tau,\omega, (\tilde{u}_\tau,\tilde{v}_\tau))=(\tilde{u}(t+\tau,\tau,\vartheta_{-\tau}\omega,\tilde{u}_\tau),\tilde{v}(t+\tau,\tau,\vartheta_{-\tau}\omega,\tilde{v}_\tau))\notag\\
 &=(z^{-1}(t+\tau,\vartheta_{-\tau}\omega)u(t+\tau,\tau,\vartheta_{-\tau}\omega,  u_\tau),z^{-1}(t+\tau,\vartheta_{-\tau}\omega)v(t+\tau,\tau,\vartheta_{-\tau}\omega,  v_\tau) ),
\end{align} where $u_\tau=z(\tau,\omega)\tilde{u}_\tau$ and $v_\tau=z(\tau,\omega)\tilde{v}_\tau$.

Suppose that for every $\tau\in\mathbb{R}$ and $\omega\in \Omega$.
\begin{align} \label{D}
\lim\limits_{t\rightarrow+\infty}e^{-\delta_1 t}\|D(\tau-t,\vartheta_{-t}\omega)\|_{L^2(\mathbb{R}^N)\times L^2(\mathbb{R}^N)}^2=0,
\end{align} where $0<\delta_0<\delta_1< \delta=\min\{\lambda,\sigma\}$.  Denote by $\mathcal{D}_{\delta}$ the collection of all families of nonempty subsets of $L^2(\mathbb{R}^N)\times L^2(\mathbb{R}^N)$ such that (\ref{D}) holds. Then it is obvious that $\mathcal{D}_{\delta}$ is inclusion closed.

 We emphasize that the choices of  the constants $\delta_1,\delta_0$ in (\ref{D}) and (\ref{a5}) respectively are different from the ones used in \cite{Adili1}. It makes us
 omit the additional assumption
\begin{align}\label{Assum}
 \lim_{t\rightarrow-\infty}e^{\delta_1}\int\limits_{-\infty}^0e^{\delta s}(\|g(s,.)\|_{L^2(\mathbb{R}^N)}^2+\|h(s,.)\|_{L^2(\mathbb{R}^N)}^2)ds=0,
\end{align}
which is intrinsically used in \cite{Adili1}, see the detailed proof of Lemma 4.1 in the following section.

Note that Adili and Wang \cite{Adili1} established the existence and upper semi-continuous of pullback attractors for problem (\ref{FN1})-(\ref{FN3}) in
$L^{2}(\mathbb{R}^N)\times L^2(\mathbb{R}^N)$. In this paper, we obtain an identical result in $L^p(\mathbb{R}^N)\times L^2(\mathbb{R}^N)$, but we do not increase the restrictions (except that $\psi_1\in L^{p/2}(\mathbb{R}^N)$ as in (\ref{a1}))
 on the nonlinearity $f$. On the contrary,  the restrictive assumption (\ref{Assum}) on the non-autonomous terms  $g$ and $h$ given in \cite{Adili1} is
 omitted.  Furthermore, we construct a unique random equilibrium for this system when  some additional assumptions on the
physical parameters are added.

\section{Existence  and upper semi-continuous of attractors in $L^{p}\times L^{2}$}

From now on, we assume  without loss of generality that $\varepsilon\in(0,a]$ for any $a>0$.
Consider that $e^{-a|\omega(s)|}\leq z_\varepsilon(s,\omega)=e^{-\varepsilon \omega(s)}\leq e^{a|\omega(s)|}$
 for $\varepsilon\in I$, and $\omega(.)$ is continuous  on $[-2,0]$. Then there exist two positive random constants
 $E=E(\omega)$ and  $F=F(\omega)$ such that for each $\omega\in\Omega$,
\begin{align} \label{EE}
&E\leq z_\varepsilon(s,\omega)\leq F\ \ \ \  \mbox{for all}\ s\in[-2,0]\ \mbox{and}\ \varepsilon\in (0,a].
\end{align}
Hereafter, we denote by $\|.\|$ and $ \|.\|_p$ the norms in $L^2(\mathbb{R}^N)$ and $L^p(\mathbb{R}^N) (p>2)$, respectively.
Throughout this paper, the number $c$  is a generic positive constant independent of $\tau,\omega,D$ and $\varepsilon$ in any place, which may vary its
values everywhere.

\subsection{Uniform absorption and uniformly asymptotic compactness in $L^2\times L^2$}

This subsection is concern with some uniform estimates of solutions on a certain compact interval $[\tau-1,\tau]$ for $\tau\in\mathbb{R}$. The uniform absorption of the family
of random cocycles $\varphi_{\varepsilon}$ is proved.  Note that the notations $(u,v)$, $(\tilde{u},\tilde{v})$, and $\varphi$ are the abbreviations of
$(u_\varepsilon,v_\varepsilon)$, $(\tilde{u}_\varepsilon,\tilde{v}_\varepsilon)$ and $\varphi_\varepsilon$ respectively, where the later implies the dependence of solutions on $\varepsilon$, omitting the subscript $\varepsilon$ if there is no confusion.\\

\textbf{Lemma 4.1.}
 \emph{Assume that (\ref{a1})-(\ref{a5}) holds and $a>0$. Given  $\tau\in\mathbb{R}, \omega\in\Omega$ and  $D=\{D(\tau,\omega);\tau\in\mathbb{R},\omega\in\Omega\}\in\mathcal{D}_\delta$, then there exists a constant  $T=T(\tau,\omega,D)\geq2$
 such that for all $t\geq T$, $\varepsilon\in (0,a]$, and $(\tilde{u}_{\tau-t},\tilde{v}_{\tau-t})\in D(\tau-t,\vartheta_{-t}\omega)$,  the solution $(u_\varepsilon,v_\varepsilon)$ of problem (\ref{pr1})-(\ref{pr2})  satisfies
\begin{align} \label{Borrow1}
\|(u_\varepsilon(\xi,\tau-t,\vartheta_{-\tau}\omega, {u}_{\tau-t}),v_\varepsilon(\xi,\tau-t,\vartheta_{-\tau}\omega,{v}_{\tau-t}))\|^2
\leq ce^{2\varepsilon\omega(-\tau)}(1+L_\varepsilon(\tau,\omega)),\xi\in [\tau-1,\tau],
\end{align}
\begin{align} \label{Borrow11}
&\|(\tilde{u}_\varepsilon(\xi,\tau-t,\vartheta_{-\tau}\omega,\tilde{u}_{\tau-t}), \tilde{v}_\varepsilon(\xi,\tau-t,\vartheta_{-\tau}\omega,\tilde{v}_{\tau-t}))\|^2\leq c(1+L_\varepsilon(\tau,\omega)),\ \  \xi\in [\tau-1,\tau],
\end{align}
and
\begin{align} \label{Borrow2}
 \int\limits_{\tau-t}^\xi e^{\delta(s-\tau)}\Big(&\|v_\varepsilon(\tau,\tau-t,\vartheta_{-\tau}\omega, {v}_{\tau-t})\|^2\notag\\
 &+z_\varepsilon^{2-p}(s,\vartheta_{-\tau}\omega)\|u_\varepsilon(s,\tau-t,\vartheta_{-\tau}\omega, {u}_{\tau-t})\|_p^p\Big)ds\leq ce^{2\varepsilon\omega(-\tau)}(1+L_\varepsilon(\tau,\omega)),
\end{align}
where $(u_{\tau-t},v_{\tau-t})=z_\varepsilon(\tau-t,\vartheta_{-\tau}\omega)(\tilde{u}_{\tau-t},\tilde{v}_{\tau-t})$, and $L_\varepsilon(\tau,\omega)$ is given by
\begin{align} \label{Borrow22}
L_\varepsilon(\tau,\omega)=\int\limits_{-\infty}^{0}e^{\delta_{01} s+2\varepsilon|\omega(s)|}\Big(\|g(s+\tau,.)\|^2+\|h(s+\tau,.)\|^2+1\Big)ds,
\end{align}} such that $\varepsilon\rightarrow L_\varepsilon(\tau,\omega)$ is an increasing function on $(0,+\infty)$, where $0<\delta_0<\delta_{01}<\delta_1<\delta$.

\emph{In particular, the family of random cocycles
$\varphi_\varepsilon(\varepsilon\in(0,a])$ defined by (\ref{eq0}) is uniformly absorbing on $(0,a]$ for any $a>0$ in the sense of Definition 2.6.}\\

\emph{Proof}\ \  Taking the inner  products of (\ref{pr1}) and (\ref{pr2})  with $u$ and $v$, respectively,  by using (\ref{a1}), we have
\begin{align} \label{energy}
\frac{d}{dt}(\beta\|u\|^2+\alpha\|v\|^2)&+\delta(\beta\|u\|^2+\alpha\|v\|^2)\notag\\
&+2\alpha_1\beta z^{2-p}(t,\omega)\|{u}\|_p^p\leq cz^2(t,\omega)(\|g(t,.)\|^2+\|h(t,.)\|^2+\|\psi_1\|_1).
\end{align}
By applying the Gronwall lemma over the interval $[\tau-t,\xi]$ with $\xi\in[\tau-1,\tau]$ and $t\geq1$,  along with $\omega$ replaced by $\vartheta_{-\tau}\omega$, we get from (\ref{trans}) that
\begin{align} \label{un01}
&\|u(\xi, \tau-t, \vartheta_{-\tau}\omega,v_{\tau-t})\|^2+\|v(\xi, \tau-t, \vartheta_{-\tau}\omega,v_{\tau-t})\|^2\notag\\
&+\int\limits_{\tau-t}^\xi e^{-\delta(\xi-s)}(\|v(s, \tau-t, \vartheta_{-\tau}\omega,v_{\tau-t})\|^2+z^{2-p}(s,\vartheta_{-\tau}\omega)\|u(s, \tau-t, \vartheta_{-\tau}\omega,v_{\tau-t})\|_p^p)ds\notag\\
&\leq ce^{-\delta(\xi-\tau+t)}(\|u_{\tau-t}\|^2+\|v_{\tau-t}\|^2)+c\int\limits_{\tau-t}^\xi e^{-\delta(\xi-s)}z^2(s,\vartheta_{-\tau}\omega)(\|g(s,.)\|^2+\|h(s,.)\|^2+1)ds\notag\\
&\leq ce^{2\varepsilon\omega(-\tau)}\Big(e^{-\delta t}z^2(-t,\omega)(\|\tilde{u}_{\tau-t}\|^2+\|\tilde{v}_{\tau-t}\|^2)+\int\limits_{\tau-t}^\tau e^{-\delta(\tau-s)-2\varepsilon\omega(s-\tau)}(\|g(s,.)\|^2+\|h(s,.)\|^2+1)ds\Big)\notag\\
&\leq ce^{2\varepsilon\omega(-\tau)}\Big(e^{-\delta t}z^2(-t,\omega)(\|\tilde{u}_{\tau-t}\|^2+\|\tilde{v}_{\tau-t}\|^2)+\int\limits_{-t}^0 e^{\delta_{01} s-2\varepsilon\omega(s)}(\|g(s+\tau,.)\|^2+\|h(s+\tau,.)\|^2+1)ds\Big),
\end{align} where $\delta=\min\{\lambda,\sigma\}$ and $\delta_{01}<\delta$.
By (\ref{3.6}) we calculate that $\lim\limits_{t\rightarrow+\infty}e^{-(\delta-\delta_1) t}z^2(-t,\omega)=0$.  Then from the property of $\mathcal{D}_\delta$  in (\ref{D}) it follows that
\begin{align} \label{un02}
\lim_{t\rightarrow+\infty}e^{-\delta t}z^2(-t,\omega)(\|\tilde{u}_{\tau-t}\|^2+\|\tilde{v}_{\tau-t}\|^2)=0.
\end{align}
Thus  (\ref{un01}) and  (\ref{un02}) together implies  that there exists a random constant $T=T(\tau,\omega, D)\geq1$ such that
for each $\varepsilon\in(0,a]$ and all $t\geq T$,  (\ref{Borrow1}) and (\ref{Borrow2}) hold.  By (\ref{trans}) and  (\ref{Borrow1}) it is showed that (\ref{Borrow11}) hold true for all $t\geq T$.

On the other hand,  from (\ref{a5}) and (\ref{3.6}) it follow  that the integral in $L_\varepsilon(\tau,\omega)$ is meaningful and
thus $L_\varepsilon(\tau,\omega)$ is finite. Further,
\begin{align} \label{un03}
L_\varepsilon(\tau-t,\vartheta_{-t}\omega)&=\int\limits_{-\infty}^{0}e^{\delta_{01} s+2\varepsilon|\vartheta_{-t}\omega(s)|}\Big(\|g(s+\tau-t,.)\|^2+\|h(s+\tau-t,.)\|^2+1\Big)ds\notag\\
(\mbox{letting}\ s-t=s')&\leq e^{\delta_{01} t+2\varepsilon|\omega(-t)|}\int\limits_{-\infty}^{-t}e^{\delta_{01} s+2\varepsilon|\omega(s)|}\Big(\|g(s+\tau,.)\|^2+\|h(s+\tau,.)\|^2+1\Big).
\end{align}
We see from (\ref{3.6}) and the relation $\delta_{0}<\delta_{01}$ that $\lim\limits_{s\rightarrow-\infty}e^{-(\delta_{01}-\delta_0) s+2\varepsilon|\omega(s)|}=0$, so there exists a positive variable $a(\omega)$ such that
$$
0<e^{-(\delta_{01}-\delta_0) s+2\varepsilon|\omega(s)|}\leq a(\omega),\  \ s\in(-\infty,0],
$$
from which it follows that
\begin{align} \label{un04}
L_\varepsilon(\tau-t,\vartheta_{-t}\omega)&\leq a(\omega)e^{\delta_{01} t+2\varepsilon|\omega(-t)|}\int\limits_{-\infty}^{-t}e^{\delta_{0} s}\Big(\|g(s+\tau,.)\|^2+\|h(s+\tau,.)\|^2+1\Big)ds.
\end{align}
Then from (\ref{D}),(\ref{3.6}), (\ref{a5}) and (\ref{un04}) we deduce that
$$
K_\varepsilon=\{K_\varepsilon(\tau,\omega)=\{(\tilde{u},\tilde{v})\in (L^2(\mathbb{R}))^2;\  \ \|(\tilde{u},\tilde{v})\|^2\leq c(1+L_\varepsilon(\tau,\omega))\}; \tau\in\mathbb{R},\omega\in\Omega\}\in\mathcal{D}_\delta,
$$ and further the union $\overline{\cup_{\varepsilon\in (0,a]}K_\varepsilon(\tau,\omega)}\subset K_a(\tau,\omega)$. Thus
$\overline{K}=\{\overline{\cup_{\varepsilon\in (0,a]}K_\varepsilon(\tau,\omega)};\tau\in\mathbb{R},\omega\in \Omega\}\in \mathcal{D}_\delta$. The measurability of
the absorbing set $K_\varepsilon(\tau,\omega)$ follows from the measurability of the variable $L_\varepsilon(\tau,\omega)$. Finally since by (\ref{Borrow22}) and (\ref{a5}) and the relation $\delta_0<\delta_{01}$,
$$
\limsup_{\varepsilon \rightarrow0}\|K_\varepsilon(\tau,\omega)\|\leq L_0(\tau,\omega)=\int\limits_{-\infty}^{0}e^{\delta_{01} s}|\Big(\|g(s+\tau,.)\|^2+\|h(s+\tau,.)\|^2+1\Big)ds<+\infty,
$$
then we have showed the uniformly absorbing of $\varphi_\varepsilon(\varepsilon\in(0,a])$ for any $a>0.$
 $\ \ \ \ \ \ \ \ \ \ \ \ \ \ \ \ \ \  \Box$\\

The uniformly asymptotic compactness in $L^2\times L^2$ has been proved by \cite{Adili1}.\\

\textbf{Lemma 4.2.} \emph{Assume that (\ref{a1})-(\ref{a5}) hold. Then the family of random cocycles $\varphi_{\varepsilon}$ defined by
(\ref{eq0}) is uniformly pullback asymptotically compact over $\varepsilon\in(0,a]$ in $L^2(\mathbb{R}^N)\times L^2(\mathbb{R}^N)$. }

\subsection{Uniformly asymptotic compactness in $L^p\times L^2$}

In this subsection, we prove that the family of random cocycles $\varphi_\varepsilon(\varepsilon(0,a])$ is uniformly asymptotically compact in  $L^p\times L^2$.
We need to prove the $L^p$-uniform boundedness of the first component of solution $u_\varepsilon$  as well as the uniform smallness of truncation of $u_\varepsilon$ in $L^p$ norm.\\

\textbf{Lemma 4.3.} \emph{Assume that (\ref{a1})-(\ref{a5}) hold. Given  $\tau\in\mathbb{R}, \omega\in\Omega$ and  $D=\{D(\tau,\omega);\tau\in\mathbb{R},\omega\in\Omega\}\in\mathcal{D}_\delta$, then there exist some random constants $C=C(\tau,\omega)$ and $T=T(\tau,\omega,D)\geq2$
such that for all $(\tilde{u}_{\tau-t},\tilde{v}_{\tau-t})\in D(\tau-t,\vartheta_{-t}\omega)$, the solution $(u_\varepsilon,v_\varepsilon)$ of problem
 (\ref{pr1})-(\ref{pr2}) satisfies}
\begin{align} \label{}
 \sup_{t\geq T}\sup_{\xi\in[\tau-1,\tau]}\sup_{\varepsilon\in(0,a]}\|{u}_\varepsilon(\xi,\tau-t,\vartheta_{-\tau}\omega, {u}_{\tau-t})\|_p^p\leq C(\tau,\omega),
\end{align} where $C(\tau,\omega)$ is independent $\varepsilon$.\\

\emph{Proof}\ \
 Multiplying (\ref{pr1}) by $|u|^{p-2}u$ and then integrating over $\mathbb{R}^N$,  we have
\begin{align} \label{so10}
&\ \  \ \ \ \ \ \ \frac{1}{p}\frac{d}{dt}\|u\|_p^p+\lambda\|u\|_p^p\leq \alpha \int\limits_{\mathbb{R}^N}v|u|^{p-1}dx\notag\\&+z(t,\omega)\int\limits_{\mathbb{R}^N}f(x,z^{-1}(t,\omega)u)|u|^{p-2}udx+z(t,\omega)\int\limits_{\mathbb{R}^N}g(t,x)|u|^{p-2}udx.
\end{align}
From (\ref{a1}) and  $\psi\in L^{p/2}$, applying Young inequality, we obtain that
\begin{align} \label{so11}
z(t,\omega)\int\limits_{\mathbb{R}^N}f(x,z^{-1}(t,\omega)u)|u|^{p-2}udx
\leq-\alpha_1z^{2-p}(t,\omega)\|u\|_{2p-2}^{2p-2}+\frac{\lambda}{4}\|u\|_p^p+cz^{p}(t,\omega)\|\psi_1\|^{p/2}_{p/2}.
\end{align}
On the other hand,
\begin{align} \label{so12}
\alpha \int\limits_{\mathbb{R}^N}v|u|^{p-1}dx\leq \frac{1}{4}\alpha_1z^{2-p}(t,\omega)\|u\|_{2p-2}^{2p-2}+cz^{p-2}(t,\omega)\|v\|^2,
\end{align}
and
\begin{align} \label{so13}
z(t,\omega)\int\limits_{\mathbb{R}^N}g(t,x)|u|^{p-2}udx\leq \frac{1}{4}\alpha_1z^{2-p}(t,\omega)\|u\|_{2p-2}^{2p-2}+cz^{p}(t,\omega)\|g(t,.)\|^2.
\end{align}
Then combination (\ref{so10})-(\ref{so13}), it give that
\begin{align} \label{so14}
\frac{d}{dt}\|u\|_p^p+\delta\|u\|_p^p\leq cz^{p-2}(t,\omega)\|v\|^2+cz^{p}(t,\omega)(\|g(t,.)\|^2+\|\psi_1\|^{p/2}_{p/2}),
\end{align} where $\delta=\min\{\lambda,\sigma\}$. Note that $\frac{1}{\xi-\tau+2}\leq 1$ for $\xi\in[\tau-1,\tau]$.
Applying Gronwall lemma (see also Lemma 5.1 in \cite{Zhao0}) over the interval $[\tau-2,\xi]$, along with $\omega$ replaced by $\vartheta_{-\tau}\omega$, we deduce that
\begin{align} \label{so15}
\|u(\xi, \tau-t,\vartheta_{-\tau}\omega,{u}_{\tau-t})\|_p^p&\leq c\int\limits_{\tau-2}^{\tau} e^{\delta (s-\tau)} \|{u}(s, \tau-t,\vartheta_{-\tau}\omega,{u}_{\tau-t})\|_p^pds\notag\\
&+c\int\limits_{\tau-2}^{\tau} e^{\delta (s-\tau)}z^{p-2}(s,\vartheta_{-\tau})\|v(s,\tau-t,\vartheta_{-\tau}\omega,{v}_{\tau-t})\|^2ds\notag\\
&+c\int\limits_{\tau-2}^{\tau} e^{\delta (s-\tau) }z^{p}(s,\vartheta_{-\tau})(\|g(s,.)\|^2+1)ds.
\end{align}
We now  estimate every term on the right hand side of (\ref{so15}).
First from (\ref{EE}) it follows that for all $s\in [\tau-2,\tau]$ and $\varepsilon\in(0,a]$, $z^{2-p}(s,\vartheta_{-\tau}\omega)=e^{\varepsilon(2-p)\omega(-\tau)}z^{2-p}(s-\tau,\omega)\geq e^{\varepsilon(2-p)\omega(-\tau)}F^{2-p}$.
Then from (\ref{Borrow2}),  there exists $T=T(\tau,\omega,D)\geq2$, such that for all $\varepsilon\in(0,a]$ and $t\geq T$,
\begin{align} \label{so0991}
 \int\limits_{\tau-2}^\tau e^{\delta(s-\tau)}\|u(s,\tau-t,\vartheta_{-\tau}\omega, {u}_{\tau-t})\|_p^pds\leq cF^{p-2}e^{ap|\omega(-\tau)|}(1+L_a(\tau,\omega)).
\end{align}
Noticing that $z^{p-2}(s,\vartheta_{-\tau})\leq e^{(p-2)\varepsilon\omega(-\tau)}F^{p-2}$ for $s\in[\tau-2,\tau]$, then from (\ref{Borrow2}) again we see that
\begin{align} \label{so150}
\int\limits_{\tau-2}^{\tau} e^{\delta (s-\tau)}z^{p-2}(s,\vartheta_{-\tau})\|v(s,\tau-t,\vartheta_{-\tau}\omega,{v}_{\tau-t})\|^2ds\leq cF^{p-2}e^{ap|\omega(-\tau)|}(1+L_a(\tau,\omega)).
\end{align}
On the other hand, by (\ref{a5}),
\begin{align} \label{so16}
\int\limits_{\tau-2}^{\tau} e^{\delta (s-\tau) }z^{p}(s,\vartheta_{-\tau})(\|g(s,.)\|^2+1)ds\leq F^{p}e^{ap|\omega(-\tau)|}\int\limits_{-2}^{0} e^{\delta s }(\|g(s+\tau,.)\|^2+1)ds<+\infty.
\end{align}
Hence (\ref{so15})-(\ref{so16}) implies the desired.
$\ \ \ \ \ \ \Box$\\

Let $M=M(\tau,\omega)>0$. Denote by  $(u-M)_+$ the positive part of $u-M$,\emph{ i.e.},
$$ (u-M)_+=\left\{
       \begin{array}{ll}
    u-M, \ \ \mbox{if}\ u> M;\\
  0,\ \ \ \ \ \ \ \  \ \mbox{if}\ u\leq M.
       \end{array}
      \right.$$
The next lemma  will show that the unbounded part of the absolute value  $|u|$ approaches zero  in $L^p$-norm on the state domain
$\mathbb{R}^N(|u(\tau,\tau-t,\vartheta_{-\tau}\omega,u_{\tau-t})|\geq M)$ for $M$ large enough, where
\begin{align} \label{}
\mathbb{R}^N(|u(\tau,\tau-t,\vartheta_{-\tau}\omega,u_{\tau-t})|\geq M)=\{x\in \mathbb{R}^N; |u(\tau,\tau-t,\vartheta_{-\tau}\omega,u_{\tau-t})|\geq M|\}.\notag
\end{align}
Note that we need not to prove some auxiliary lemmas except Lemma 4.1 and
Lemma 4.3, see \cite{Liyangrong0,Liyangrong1,Liyangrong2,Lijia,Zhao3,Zhao4}.\\

\textbf{Lemma 4.4.} Assume that (\ref{a1})-(\ref{a5}) hold. Given  $\tau\in\mathbb{R}, \omega\in\Omega$ and  $D=\{D(\tau,\omega);\tau\in\mathbb{R},\omega\in\Omega\}\in\mathcal{D}_\delta$, then for any $\eta>0$, there exist  random constants $M=M(\tau,\omega,\eta,D)>1$  and $T=T(\tau,\omega,D)\geq2$
such that  for all $(\tilde{u}_{\tau-t},\tilde{v}_{\tau-t})\in D(\tau-t,\vartheta_{-t}\omega)$, the first component $\tilde{u}_\varepsilon$ of solutions $(\tilde{u}_\varepsilon,\tilde{v}_\varepsilon)$ of problem (\ref{FN1})-(\ref{FN3}) satisfies
\begin{align} \label{}
\sup_{t\geq T}\sup_{\varepsilon\in(0,a]}\int\limits_{\mathbb{R}^N(|\tilde{u}_\varepsilon|\geq M)}|\tilde{u}_\varepsilon(\tau,\tau-t,\vartheta_{-\tau}\omega,\tilde{u}_{\tau-t})|^p dx\leq \eta,\notag
\end{align} where $\mathbb{R}^N(|\tilde{u}_\varepsilon|\geq M)=\mathbb{R}^N(|\tilde{u}_\varepsilon(\tau,\tau-t,\vartheta_{-\tau}\omega,\tilde{u}_{\tau-t})|\geq M)$, and $M,T$ are independent of $\varepsilon$.\\

\emph{Proof} \ \  Let $s\in[\tau-1,\tau]$ and $t\geq T\geq2$, where $T$ is determined by Lemma 4.1 and Lemma 4.2.
Replacing  $\omega$ by $\vartheta_{-\tau}\omega$ in (\ref{pr1})-(\ref{pr2}), we see that
$u=u(s,\tau-t,\vartheta_{-\tau}\omega,u_{\tau-t}),v=v(s,\tau-t,\vartheta_{-\tau}\omega,v_{\tau-t})$, is a solution of the following system
\begin{align} \label{p011}
\frac{du}{ds}+\lambda u-\Delta u+\alpha v=\frac{z(s-\tau,\omega)}{z(-\tau,\omega)}f(x,\tilde{u})+\frac{z(s-\tau,\omega)}{z(-\tau,\omega)}{g}(s,x),
\end{align}
\begin{align} \label{p0112}
\frac{dv}{ds}+\sigma v-\beta u=\frac{z(s-\tau,\omega)}{z(-\tau,\omega)}{h}(s,x).
\end{align}
For fixed $\tau\in \mathbb{R}$ and $\omega\in\Omega$, we assume that $ M=M(\tau,\omega)>1$.
We  multiply (\ref{p011}) by  $(u-M)_+^{p-1}$  and integrate over $\mathbb{R}^N$ to yield
\begin{align} \label{p01}
\frac{1}{p}\frac{d}{ds}\int\limits_{\mathbb{R}^N} (u-M)_+^{p}dx&+\lambda\int\limits_{\mathbb{R}^N} u(u-M)_+^{p-1}dx-\int\limits_{\mathbb{R}^N}\Delta u(u-M)_+^{p-1}dx
\notag\\&=-\alpha\int\limits_{\mathbb{R}^N}v(u-M)_+^{p-1}dx+ \frac{z(s-\tau,\omega)}{z(-\tau,\omega)} \int\limits_{\mathbb{R}^N}f(x,\tilde{u})(u-M)_+^{p-1}dx\notag\\
&\ \ \  \ \  \ \ +\frac{z(s-\tau,\omega)}{z(-\tau,\omega)} \int\limits_{\mathbb{R}^N}{g}(s,x)(u-M)_+^{p-1}dx.
\end{align}
We now have to estimate every term in (\ref{p01}). First, it is obvious that
\begin{align} \label{p02}
-\int\limits_{\mathbb{R}^N}\Delta u(u-M)_+^{p-1}dx=(p-1)\int\limits_{\mathbb{R}^N}(u-M)_+^{p-2}|\nabla
u|^2dx\geq0,
\end{align}
\begin{align} \label{p03}
\lambda\int\limits_{\mathbb{R}^N} u(u-M)_+^{p-1}dx\geq \lambda\int\limits_{\mathbb{R}^N} (u-M)_+^{p}dx.
\end{align}
The most involved work is to calculate the nonlinearity in (\ref{p01}). Consider that for $u(s)> M$ for $\in[\tau-1,\tau]$, we have $\tilde{u}(s)=z^{-1}(s,\vartheta_{-\tau}\omega)u(s)=\frac{z(-\tau,\omega)}{z(s-\tau,\omega)}u(s)>0$, and thus
by (\ref{a1}) and (\ref{EE}), we find that for every $s\in[\tau-1,\tau]$,
\begin{align} \label{}
f(x,\tilde{u})&\leq -\alpha_1|\tilde{u}|^{p-1}+\frac{1}{\tilde{u}}\psi_1(x)\notag\\
&= -\alpha_1\Big(\frac{z(s-\tau,\omega)}{z(-\tau,\omega)}\Big)^{1-p}|u|^{p-1}+\frac{z(s-\tau,\omega)}{z(-\tau,\omega)}\frac{\psi_1(x)}{u}\notag\\
&\leq -\frac{\alpha_1z^{p-1}(-\tau,\omega)}{2F^{p-1}}M^{p-2}(u-M)- \frac{\alpha_1z^{p-1}(-\tau,\omega)}{2F^{p-1}}(u-M)^{p-1}
\notag\\
&+\frac{F}{z(-\tau,\omega)}|\psi_1(x)|(u-M)^{-1},\notag
\end{align}
 from which and (\ref{EE}) again it follows that
\begin{align} \label{p04}
&\ \ \  \ \ \ \ \ \ \ \ \ \ \ \ \ \ \  \frac{z(s-\tau,\omega)}{z(-\tau,\omega)} \int\limits_{\mathbb{R}^N}f(x,\tilde{u})(u-M)_+^{p-1}dx\notag\\&
\leq-\frac{\alpha_1Ez^{p-1}(-\tau,\omega)}{2F^{p-1}}M^{p-2}\int\limits_{\mathbb{R}^N}(u-M)_+^{p}dx-
\frac{\alpha_1Ez^{p-1}(-\tau,\omega)}{2F^{p-1}}\int\limits_{\mathbb{R}^N}(u-M)_+^{2p-2}dx\notag\\
&\ \  \ \ \ \ \ \ \ \ \ \ \ \ \ \ \ \ \ +\frac{F^2}{z^2(-\tau,\omega)}\int\limits_{\mathbb{R}^N}|\psi_1(x)|(u-M)_+^{p-2}dx\notag\\
&\leq-\frac{\alpha_1Ez^{p-1}(-\tau,\omega)}{2F^{p-1}}M^{p-2}\int\limits_{\mathbb{R}^N}(u-M)_+^{p}dx-\frac{\alpha_1Ez^{p-1}(-\tau,\omega)}{2F^{p-1}}\int\limits_{\mathbb{R}^N}(u-M)_+^{2p-2}dx\notag\\
& \ \ \ \ \  \ \ \ \ \ \ \ \  \ \ \ \ +\frac{1}{2}\lambda\int\limits_{\mathbb{R}^N}(u-M)_+^{p}dx+\frac{cF^p}{z^p(-\tau,\omega)}\int\limits_{\mathbb{R}^N(u\geq M)}|\psi_1(x)|^{p/2}dx,
\end{align} in which  we have used the Young inequality in the last term, and here $c=(\frac{2}{\lambda})^{p-2/2}$.
On the other hand by using the Young inequality again, we get that for $s\in [\tau-1,\tau]$,
\begin{align} \label{p05}
\Big|\frac{z(s-\tau,\omega)}{z(-\tau,\omega)}& \int\limits_{\mathbb{R}^N}{g}(s,x)(u-M)_+^{p-1}dx\Big|\leq\frac{F}{z(-\tau,\omega)}\Big| \int\limits_{\mathbb{R}^N}{g}(s,x)(u(s)-M)_+^{p-1}dx\Big|\notag\\
&\leq \frac{\alpha_1Ez^{p-1}(-\tau,\omega)}{4F^{p-1}}\int\limits_{\mathbb{R}^N}(u-M)_+^{2p-2}dx+\frac{F^{p+1}}{\alpha_1Ez^{p+1}(-\tau,\omega)}\int\limits_{\mathbb{R}^N(u(s)\geq M)}{g}^2(s,x)dx,
\end{align}
and
\begin{align} \label{p051}
\Big|-\alpha\int\limits_{\mathbb{R}^N}v(u-M)_+^{p-1}dx\Big|&\leq\frac{\alpha_1Ez^{p-1}(-\tau,\omega)}{4F^{p-1}}\int\limits_{\mathbb{R}^N}(u-M)_+^{2p-2}dx\notag\\
&+\frac{\alpha^2F^{p-1}}{\alpha_1Ez^{p-1}(-\tau,\omega)}
\int\limits_{\mathbb{R}^N(u(s)\geq M)}v^2dx.
\end{align}
For convenience of calculations, we introduce the following notations:
\begin{align} \label{p0510}
k=k(\tau,\omega,M)=\frac{\alpha_1Ee^{-(p-1)a|\omega(-\tau)|}}{2F^{p-1}}M^{p-2};
\end{align} which is increasing to infinite in $M$ for $p>2$,
and
\begin{align} \label{p05101}
G(\tau,\omega)
=\max\Big\{\frac{F^{p+1}e^{a(p+1)|\omega(-\tau)|}}{\alpha_1E};\frac{\alpha^2F^{p-1}e^{a(p-1)|\omega(-\tau)|}}{\alpha_1E};(\frac{2}{\lambda})^{\frac{p-2}{2}}F^pz^{ap|\omega(-\tau)|}\Big\},
\end{align} which is a nonnegative random constant depending only on $\tau,\omega$.
Combination (\ref{p01})-(\ref{p051}) and using the notations (\ref{p0510})-(\ref{p05101}), we obtain that
\begin{align} \label{p06}
\frac{d}{ds}\int\limits_{\mathbb{R}^N} (u(s)-M)_+^{p}dx&+
k\int\limits_{\mathbb{R}^N}(u(s)-M)_+^{p}dx\leq G(\tau,\omega)(\|{g}(s,.)\|^2+\|v\|^2+1),
\end{align} where $s\in[\tau-1,\tau]$.
Applying Gronwall lemm (also see Lemma 5.1 in \cite{Zhao0} ) over $[\tau-1,\tau]$, by Lemma 4.1 and Lemma 4.3, we
find that for all $t\geq T$ and $\varepsilon\in(0,a]$,
 \begin{align} \label{p07}
\int\limits_{\mathbb{R}^N} \Big(u(\tau,\tau-t,&\vartheta_{-\tau}\omega,u_{\tau-t})-M\Big)_+^{p}dx
\leq \int\limits_{\tau-1}^\tau e^{k(s-\tau)}\|u(s,\tau-t,\vartheta_{-\tau}\omega,u_0)\|_p^pds\notag\\
&+G(\tau,\omega)\Big(\int\limits_{\tau-1}^\tau e^{k(s-\tau)} \|v(s,\tau-t,\vartheta_{-\tau}\omega,v_{\tau-t})\|^2ds+\int\limits_{\tau-1}^\tau e^{k(s-\tau)} (\|g(s,.)\|^2+1)ds\Big)\notag\\
&\leq \frac{C(\tau,\omega)+cG(\tau,\omega)e^{2a|\omega(-\tau)|}(1+L_a(\tau,\omega))+G(\tau,\omega)}{k}\notag\\
&+G(\tau,\omega)\int\limits_{\tau-1}^\tau e^{k(s-\tau)}\|g(s,.)\|^2ds.
\end{align}
For fixed $\tau\in\mathbb{R}$ and $\omega\in\Omega$,  the first term in the last inequality of (\ref{p07}) varies only with the number $k$, but $k$ is a large number
as $M\rightarrow+\infty$. Therefore, we get that this term converges to zero when $M$ goes to infinite. It remains to prove that the second term vanishes for $M$ large enough.
 First, choosing a large $M$ such that $k=k(\tau,\omega,M)>\delta_{0}$ ($\delta_{0}$ is as in (\ref{a5})) and taking $\varsigma\in(0,1)$, we have
\begin{align} \label{p08}
\int\limits_{\tau-1}^\tau e^{k(s-\tau)}\|g(s,.)\|^2ds&=\int\limits_{\tau-1}^{\tau-\varsigma} e^{k(s-\tau)}\|g(s,.)\|^2ds+\int\limits_{\tau-\varsigma}^{\tau} e^{k(s-\tau)}\|g(s,.)\|^2ds\notag\\
&=e^{-k\tau}\int\limits_{\tau-1}^{\tau-\varsigma} e^{(k-\delta_{0})s}e^{\delta s}\|g(s,.)\|^2ds+e^{-k\tau}\int\limits_{\tau-\varsigma}^{\tau} e^{ks}\|g(s,.)\|^2ds\notag\\
&\leq e^{-k\varsigma} e^{\delta_{0}(\varsigma-\tau)}\int\limits_{-\infty}^{\tau} e^{\delta_{0} s}\|g(s,.)\|^2ds+\int\limits_{\tau-\varsigma}^{\tau}\|g(s,.)\|^2ds.
\end{align}
By (\ref{a5}), the first term above vanishes as $k\rightarrow+\infty$, and by $g\in L^2_{loc}(\mathbb{R}, L^2(\mathbb{R}^N))$  we can choose $\varsigma$ small enough such that the second term in (\ref{p08}) is small.
In terms of these arguments, from (\ref{p07})and (\ref{p08}) we have proved that
\begin{align} \label{p077}
\sup_{t\geq T}\sup_{\varepsilon\in(0,a]}\int\limits_{\mathbb{R}^N} \Big(u(\tau,\tau-t,\vartheta_{-\tau}\omega,u_{\tau-t})-M\Big)_+^{p}dx\rightarrow0,
\end{align} when $M\rightarrow+\infty$.
Therefore, for any $\eta>0$, there exists  $M_1=M_1(\tau,\omega,\eta,D)>1$ large enough such that
\begin{align} \label{p10}
\sup_{t\geq T}\sup_{\varepsilon\in(0,a]}\int\limits_{\mathbb{R}^N} \Big(u(\tau,\tau-t,\vartheta_{-\tau}\omega,u_{\tau-t})-M_1\Big)_+^{p}dx&\leq e^{-ap|\omega(-\tau)|}\frac{\eta}{2^{p+1}}.
\end{align}
If $u(\tau,\tau-t,\vartheta_{-\tau}\omega,u_{\tau-t})\geq 2M_1$, then $u(\tau,\tau-t,\vartheta_{-\tau}\omega,u_{\tau-t})-M_1\geq \frac{u(\tau,\tau-t,\vartheta_{-\tau}\omega,u_{\tau-t})}{2},$ so
by (\ref{p10}) it infer us  that
\begin{align} \label{p11}
&\ \  \ \sup_{t\geq T}\sup_{\varepsilon\in(0,a]}\int\limits_{\mathbb{R}^N(u(\tau)\geq 2M_1)}
|{u}(\tau,\tau-t,\vartheta_{-\tau}\omega,{u}_{\tau-t})|^{p}dx\leq  e^{-ap|\omega(-\tau)|}\frac{\eta}{2}.
\end{align}
We see from (\ref{trans}) that $\tilde{u}(\tau,\tau-t,\vartheta_{-\tau}\omega,\tilde{u}_{\tau-t})=z(-\tau,\omega)u(\tau,\tau-t,\vartheta_{-\tau}\omega,u_{\tau-t})$. Then in terms of the fact that $e^{-a|\omega(-\tau)|}\leq z(-\tau,\omega)=e^{-\varepsilon\omega(-\tau)}\leq e^{a|\omega(-\tau)|}$ for all $\varepsilon\in(0,a]$, it induce that
$\mathbb{R}^N(\tilde{u}(\tau,\tau-t,\vartheta_{-\tau}\omega,\tilde{u}_{\tau-t})\geq 2M_1e^{a|\omega(-\tau)|})\subseteq \mathbb{R}^N(u(\tau,\tau-t,\vartheta_{-\tau}\omega,u_{\tau-t})\geq 2M_1).$
This along with  (\ref{p11}) implies that,
\begin{align} \label{p80}
 &\sup_{t\geq T}\sup_{\varepsilon\in(0,a]}\int\limits_{\mathbb{R}^N(\tilde{u}(\tau)\geq 2M_1e^{a|\omega(-\tau)|})}|\tilde{u}(\tau,\tau-t,\vartheta_{-\tau}\omega,\tilde{u}_{\tau-t})|^{p}dx\notag\\
&\ \ \ \ \ \ \leq\sup_{t\geq T}\sup_{\varepsilon\in(0,a]} e^{ap|\omega(-\tau)|}\int\limits_{\mathbb{R}^N(u(\tau)\geq 2M_1)}|u(\tau,\tau-t,\vartheta_{-\tau}\omega,u_{\tau-t})|^{p}dx
\leq \frac{\eta}{2}.
\end{align}
Similarly, we can deduce that there exists  $M_2=M_2(\tau,\omega,\eta,D)>0$ large enough such that
\begin{align} \label{p81}
\sup_{t\geq T}\sup_{\varepsilon\in(0,1]}\int\limits_{\mathbb{R}^N(\tilde{u}(\tau)\leq - 2M_2e^{a|\omega(-\tau)|})}&|\tilde{u}(\tau,\tau-t,\vartheta_{-\tau}\omega,\tilde{u}_{\tau-t})|^{p}dx\leq \frac{\eta}{2}.
\end{align} Put $M=\max\{M_1,M_2\}\times e^{a|\omega(-\tau)|} $. Then (\ref{p80}) and (\ref{p81}) together imply the desired.
 $\ \ \  \ \ \ \ \ \ \Box$\\

\textbf{Lemma 4.5.} \emph{Assume that (\ref{a1})-(\ref{a5}) hold. Then  for every  $\tau\in\mathbb{R}, \omega\in\Omega$,
 $\{\tilde{u}_{\varepsilon_n}(\tau, \tau-t_n,\vartheta_{-t_n}\omega, \tilde{u}_{0,n})\}$ has a convergent subsequence in
 $L^p(\mathbb{R}^N)$ whenever $\varepsilon_n\in (0,a]$, $t_n\rightarrow+\infty$ and $(\tilde{u}_{0,n},\tilde{v}_{0,n})\in D(\tau-t_n,\vartheta_{-t_n}\omega)\in \mathcal{D}_{\delta}$.}\\

\emph{Proof.} \ \ Denote by $\tilde{u}_n(\tau)=\tilde{u}_{\varepsilon_n}(\tau,\tau-t_n,\vartheta_{-\tau}\omega,\tilde{u}_{0,n})$.
 From Lemma 4.3,  for any $\eta>0$, there exist random constants $M=M(\tau,\omega,\eta,D)>1$  and $\mathcal{Z}_1=\mathcal{Z}_1(\tau,\omega,D)\in \mathbb{Z}^+$
such that the solution $\tilde{u}_n(\tau)$  satisfies that for all $n\geq \mathcal{Z}_1$,
\begin{align} \label{p20}
\sup_{\varepsilon_n\in(0,a]}\int\limits_{\mathbb{R}^N(|\tilde{u}_n(\tau)|\geq M)}|\tilde{u}_n(\tau)|^p dx\leq \frac{\eta^p}{2^{p+2}}.
\end{align}
On the other hand, Lemma 4.2 also implies that
there exists a $\mathcal{Z}_2=\mathcal{Z}_1(\tau,\omega,B)\in \mathbb{Z}^+$  such that for all $n,m\geq \mathcal{Z}_2$,
\begin{align} \label{p21}
\int\limits_{\mathbb{R}^N}|\tilde{u}_n(\tau)-\tilde{u}_m(\tau)|^2dx\leq \frac{1}{(2M)^{p-2}}\frac{\eta^p}{4},
\end{align} whenever $\varepsilon_n,\varepsilon_m\in(0,a]$. Here $M$ is as in (\ref{p20}). We then decompose the entire space $\mathbb{R}^N$  by $\mathbb{R}^N=\mathcal{O}_1\cup \mathcal{O}_2\cup
\mathcal{O}_3\cup \mathcal{O}_4$, where
\begin{align}
&\mathcal{O}_1=\mathbb{R}^N (|\tilde{u}_n(\tau)|\leq M)\cap \mathcal{O}^N( |\tilde{u}_m(\tau)|\leq
M);\ \ \mathcal{O}_2=\mathbb{R}^N (|u_n(\tau)|\geq M)\cap
\mathbb{R}^N(|\tilde{u}_m(\tau)|\leq
M);\notag\\
&\mathcal{O}_3=\mathbb{R}^N (|\tilde{u}_n(\tau)|\leq M)\cap \mathbb{R}^N( |\tilde{u}_m(\tau)|\geq
M);\ \ \mathcal{O}_4=\mathbb{R}^N (|\tilde{u}_n(\tau)|\geq M)\cap
\mathbb{R}^N(|\tilde{u}_m(\tau)|\geq M).\notag
 \end{align}
We now put $\mathcal{Z}=\max\{\mathcal{Z}_1,\mathcal{Z}_2\}$. Then for all $n,m\geq \mathcal{Z}$, (\ref{p20}) and (\ref{p21}) hold true.
By (\ref{p21}),  we have
\begin{align}\label{p22}
 \int\limits_{\mathcal{O}_1}|\tilde{u}_n(\tau)-\tilde{u}_m(\tau)|^pdx
&\leq \int\limits_{\mathbb{R}^N(|\tilde{u}_n(\tau)-\tilde{u}_m(\tau)|\leq
2M)}|\tilde{u}_n(\tau)-\tilde{u}_m(\tau)|^pdx
\notag\\&\leq (2M)^{p-2}\|\tilde{u}_n(\tau)-\tilde{u}_m(\tau)\|^2
\leq (2M)^{p-2}.(2M)^{2-p}(\frac{\eta^p}{4})=\frac{\eta^p}{4}.
 \end{align}
On the other hand, according to (\ref{p20}),
\begin{align}\label{p23}
\int\limits_{\mathcal{O}_2}|\tilde{u}_n(\tau)-\tilde{u}_m(\tau)|^pdx &\leq 2^p\int\limits_{\mathbb{R}^N(|\tilde{u}_n(\tau)|\geq
M)}|\tilde{u}_n(\tau)|^pdx \leq\frac{\eta^p}{4},
 \end{align}
\begin{align}\label{p24}
\int\limits_{\mathcal{O}_3}|\tilde{u}_n(\tau)-\tilde{u}_m(\tau)|^pdx &\leq 2^p\int\limits_{\mathbb{R}^N(|\tilde{u}_m(\tau)|\geq
M)}|\tilde{u}_m(\tau)|^pdx \leq\frac{\eta^p}{4},
 \end{align}
\begin{align}\label{p25}
\int\limits_{\mathcal{O}_4}|\tilde{u}_n(\tau)-\tilde{u}_m(\tau)|^pdx &\leq
2^{p-1}\Big(\int\limits_{\mathbb{R}^N(|u(\tau)|\geq
M)}|\tilde{u}_n(\tau)|^pdx+\int\limits_{\mathbb{R}^N(|\tilde{u}_m(\tau)|\geq M)}|\tilde{u}_m(\tau)|^pdx\Big)
\leq\frac{\eta^p}{4}.
 \end{align}
It follows form (\ref{p22})-(\ref{p25}) that
$$\|\tilde{u}_n(\tau)-\tilde{u}_m(\tau)\|_p\leq\eta\ \ \ \mbox{for all}\ n,m\geq\mathcal{Z},$$
whenever $\varepsilon_n,\varepsilon_m\in(0,a]$, which shows that $\{\tilde{u}_n(\tau)\}$ also has a convergent subsequence  in $L^p(\mathbb{R}^N)$. Then the proof is concluded.
 $\ \ \  \ \ \ \ \ \ \Box$\\

By Lemma 4.2 and Lemma 4.5 we immediately have
\\

\textbf{Lemma  4.6.} \emph{Assume that (\ref{a1})-(\ref{a5}) hold. Then the family of random cocycles $\varphi_{\varepsilon}$ defined by
(\ref{eq0}) is uniformly pullback asymptotically compact over $\varepsilon\in(0,a]$ in $L^p(\mathbb{R}^N)\times L^2(\mathbb{R}^N)$. In particular, each $\varphi_\varepsilon$
has a unique  $(L^2(\mathbb{R}^N)\times L^2(\mathbb{R}^N),L^p(\mathbb{R}^N)\times L^2(\mathbb{R}^N))$-pullback attractor $\mathcal{A}_\varepsilon$}\\

\emph{Proof}\ \ The uniformly pullback asymptotic compactness is followed from  Lemma 4.2 and Lemma 4.4, and the existence and uniqueness of bi-spatial pullback attractor are from Theorem 2.9, or Theorem 2.10.
$\ \ \  \ \ \ \ \ \ \Box$\\

\subsection{Convergence of the family $\varphi_\varepsilon$ on $(0,a]$ in $L^{2}\times L^2$}

This subsection deal with the convergence of solutions at any intension $\varepsilon$ of noise. The convergence at zero has been shown
 by \cite{Adili1}. Here we need to prove it also converges at any $\varepsilon>0$. To this end,
 the following assumption on the nonlinearity $f$   as in \cite{Adili2} is also required. That is, for all $x\in\mathbb{R}^N$ and
$s\in\mathbb{R}$,
\begin{align}\label{a6}
 \Big|\frac{\partial}{\partial s}(x,s)\Big|\leq \alpha_4|s|^{p-2}+\psi_4(x),
\end{align}
where $\alpha_4>0$,  $\psi_4\in L^\infty(\mathbb{R}^N)$ if $p=2$ and $\psi_4\in L^{\frac{p}{p-2}}(\mathbb{R}^N)$ if $p>2$. We  need further to assume that $\psi_2\in L^q(\mathbb{R}^N)$, where $\psi_2$ is as in (\ref{a2}) and $q=\frac{p}{p-1}$ is conjugation of $p$.\\

To begin with, from (\ref{energy}), it  is very easy to derive the following inequality. \\

\textbf{ Lemma  4.7.} \emph{ Assume that (\ref{a1})-(\ref{a5}) hold.
 Then for each $\tau\in\mathbb{R}$,  $\omega\in \Omega$ and  $(\tilde{u}_\varepsilon(\tau),\tilde{v}_\varepsilon(\tau))\in L^2(\mathbb{R}^N)\times L^2(\mathbb{R}^N)$, the solution $({u},{v})$ of problem (\ref{pr1}) satisfies for all $t\geq \tau$,
 \begin{align}
 &\|u_\varepsilon(t,\tau,\omega,u(\tau))\|^2+ \|v_\varepsilon(t,\tau,\omega,v(\tau))\|^2\notag\\
 &\ \  \ \ \ \  \ \ \ \ \  +\int\limits_{\tau}^t e^{\delta(s-t)}\Big(\|v_\varepsilon(s,\tau,\omega,v_\varepsilon(\tau))\|^2+z_\varepsilon^2(s,\omega)\|\tilde{u}_\varepsilon(s,\tau,\omega,\tilde{u}_\varepsilon(\tau))\|^p_p\Big)ds\notag\\
 &\leq z_\varepsilon^2(\tau,\omega)(\|\tilde{u}_\varepsilon(\tau)\|^2+\|\tilde{v}_\varepsilon(\tau)\|^2)+c\int\limits_{\tau}^tz_\varepsilon^2(s,\omega)(\|g(s,.)\|^2+\|h(s,.)\|^2+1)ds.\notag
\end{align}}

Then by applying Lemma 4.7 we have\\

\textbf{ Lemma  4.8.}  \emph{ Assume that (\ref{a1})-(\ref{a5}) and (\ref{a6}) hold. Let $(\tilde{u}_\varepsilon,\tilde{v}_\varepsilon)$ be the solution of
 problem (\ref{pr1})-(\ref{pr2}) with initial data $(\tilde{u}_{\varepsilon,\tau},\tilde{v}_{\varepsilon,\tau})$. Assume that $\varepsilon\rightarrow\varepsilon_{0}$ and
 $\|(\tilde{u}_{\varepsilon,\tau},\tilde{v}_{\varepsilon,\tau})-(\tilde{u}_{\varepsilon_{0},\tau},\tilde{v}_{\varepsilon_0,\tau})\|\rightarrow0$  for $\varepsilon,\varepsilon_{0}\in(0,a]$.
 Then for each $\tau\in\mathbb{R}$,  $\omega\in \Omega, T>0$ and every $t\in [\tau,\tau+T]$,
\begin{align}\label{sss1}
 \lim\limits_{\varepsilon\rightarrow\varepsilon_{0}}\|(\tilde{u}_{\varepsilon}(t,\tau,\omega,\tilde{u}_{\varepsilon,\tau}),\tilde{v}_{\varepsilon}(t, \tau,\omega,\tilde{v}_{\varepsilon,\tau}))-(\tilde{u}_{\varepsilon_{0}}(t, \tau, \omega, \tilde{u}_{\varepsilon_0,\tau}),\tilde{v}_{\varepsilon_{0}}(t,\tau,\omega,\tilde{v}_{\varepsilon_0,\tau}))\|=0.
 \end{align}
 In particular,  let $(\tilde{u},\tilde{v})$ be the solution of problem (\ref{pr1})-(\ref{pr2}) for $\varepsilon=0$ with initial data $(\tilde{u}_\tau,\tilde{v}_\tau)$.
 Assume that $\varepsilon\rightarrow{0}$ and
 $\|(\tilde{u}_{\varepsilon,\tau},\tilde{v}_{\varepsilon,\tau}))-(\tilde{u}_{\tau},\tilde{v}_{\tau})\|\rightarrow0$. Then for each $\tau\in\mathbb{R}$,  $\omega\in \Omega, T>0$ \\
\begin{align}\label{sss2}
 \lim\limits_{\varepsilon\rightarrow{0}}\|(\tilde{u}_{\varepsilon}(t,\tau,\omega,\tilde{u}_{\varepsilon,\tau}),\tilde{v}_{\varepsilon}(t, \tau,\omega,\tilde{v}_{\varepsilon,\tau}))-(\tilde{u}(t, \tau, \tilde{u}_{\tau}),\tilde{v}(t,\tau,\tilde{v}_{\tau}))\|=0.
 \end{align}}

\emph{Proof}\ \  Put $U=U(t)=u_{\varepsilon}(t,\tau,\omega,u_{\varepsilon,\tau})-u_{\varepsilon_{0}}(t,\tau,\omega,u_{\varepsilon_0,\tau})$ and $V=V(t)=v_{\varepsilon}(t,\tau,\omega,v_{\varepsilon,\tau})-v_{\varepsilon_{0}}(t,\tau,\omega,v_{\varepsilon_0,\tau})$. Then we get the following system:
\begin{equation}
\begin{cases}
\frac{dU}{dt}+\lambda U-\Delta U+\alpha V=e^{-\varepsilon\omega(t)}f(x,e^{\varepsilon\omega(t)}u_\varepsilon)-e^{-\varepsilon_{0}\omega(t)}f(x,e^{\varepsilon_{0}\omega(t)}u_{\varepsilon_{0}})\\
  \ \ \ \ \ \ \ \ \ \ \ \ \ \ \  \ \ \ \ \ \  \ \ \ \ \ \ \ \ +(e^{-\varepsilon\omega(t)}-
e^{-\varepsilon_{0}\omega(t)})g(t,x),\\
\frac{dV}{dt}+\sigma V-\beta U=(e^{-\varepsilon\omega(t)}-e^{-\varepsilon_{0}\omega(t)})h(t,x),
\end{cases}
\end{equation}
where $u_\varepsilon=u_\varepsilon(t)=u_{\varepsilon}(t,\tau,\omega,u_{\varepsilon,\tau})$.
Let $\eta$ be a small positive number. Since $\omega$ is continuous on $\mathbb{R}$, then there exists a $\chi=\chi(\tau,\omega,\eta,T)>0$ such that
for every $\varepsilon\in (\varepsilon_{0}-\chi,\varepsilon_{0}+\chi)\subset(0,a]$ and $t\in [\tau,\tau+T]$,
\begin{align}\label{U02}
|e^{\varepsilon\omega(t)}-e^{\varepsilon_{0}\omega(t)}|+|e^{-\varepsilon\omega(t)}-e^{-\varepsilon_{0}\omega(t)}|\leq \eta.
\end{align}
By (4.48), we deduce that
\begin{align}\label{U03}
&\frac{1}{2}\frac{d}{dt}(\beta\|U\|^2+\alpha\|V\|^2)+\lambda \beta\|U\|^2+\sigma\alpha\|V\|^2\notag\\
&\leq \int\limits_{\mathbb{R}^N}\Big(e^{-\varepsilon\omega(t)}f(x,e^{\varepsilon\omega(t)}u_\varepsilon)-e^{-\varepsilon_{0}\omega(t)}f(x,e^{\varepsilon_{0}\omega(t)}u_{\varepsilon_{0}})\Big)Udx
\notag\\
&+(e^{-\varepsilon\omega(t)}-
e^{-\varepsilon_{0}\omega(t)})\int\limits_{\mathbb{R}^N} g(t,x)Udx +(e^{-\varepsilon\omega(t)}-
e^{-\varepsilon_{0}\omega(t)})\int\limits_{\mathbb{R}^N} h(t,x)Vdx.
\end{align}
 The first term on the right  hand side of  (\ref{U03}) is  rewritten as
\begin{align}\label{U04}
& \int\limits_{\mathbb{R}^N}\Big(e^{-\varepsilon\omega(t)}f(x,e^{\varepsilon\omega(t)}u_\varepsilon)-e^{-\varepsilon_{0}\omega(t)}f(x,e^{\varepsilon_{0}\omega(t)}u_{\varepsilon_{0}})\Big)Udx\notag\\
&=e^{-\varepsilon\omega(t)}\int\limits_{\mathbb{R}^N}\Big(f(x,e^{\varepsilon\omega(t)}u_\varepsilon)-f(x,e^{\varepsilon_{0}\omega(t)}u_{\varepsilon_{0}})\Big)Udx
+(e^{-\varepsilon\omega(t)}-e^{-\varepsilon_{0}\omega(t)})\int\limits_{\mathbb{R}^N}f(x,e^{\varepsilon_{0}\omega(t)}u_{\varepsilon_{0}})Udx\notag\\
&\ \  \ \ \ \   =e^{-\varepsilon\omega(t)}\int\limits_{\mathbb{R}^N}\frac{\partial }{\partial s}f(x,s)(e^{\varepsilon\omega(t)}u_\varepsilon-e^{\varepsilon_{0}\omega(t)}u_{\varepsilon_{0}})Udx
+(e^{-\varepsilon\omega(t)}-e^{-\varepsilon_{0}\omega(t)})\int\limits_{\mathbb{R}^N}f(x,e^{\varepsilon_{0}\omega(t)}u_{\varepsilon_{0}})Udx\notag\\
&\ \  \ \ \ \  \ \ \ \ \ \ \ \ \ \ =\int\limits_{\mathbb{R}^N}\frac{\partial }{\partial s}f(x,s)U^2dx+(e^{-\varepsilon_{0}\omega(t)}-e^{-\varepsilon\omega(t)})\int\limits_{\mathbb{R}^N}
\frac{\partial }{\partial s}f(x,s)\tilde{u}_{\varepsilon_{0}}Udx\notag\\
& \ \ \ \ \ \ \ \ \ \ \ \ \  \ \ \ \ \ \ \ \ \ \ \ \ \ \ \ +(e^{-\varepsilon\omega(t)}-e^{-\varepsilon_{0}\omega(t)})\int\limits_{\mathbb{R}^N}f(x,e^{\varepsilon_{0}\omega(t)}u_{\varepsilon_{0}})Udx.
\end{align}
By  (\ref{a6}) and (\ref{U02}), the second term on the right hand side of ({\ref{U04}}) is bounded by
\begin{align}\label{U05}
&(e^{-\varepsilon_{0}\omega(t)}-e^{-\varepsilon\omega(t)})\int\limits_{\mathbb{R}^N}
\frac{\partial }{\partial s}f(x,s)\tilde{u}_{\varepsilon_{0}}Udx\notag\\
& \ \ \  \ \ \  \ \ \ \ \leq |e^{-\varepsilon_{0}\omega(t)}-e^{-\varepsilon\omega(t)}|\int\limits_{\mathbb{R}^N}\Big(\alpha_4|\tilde{u}_\varepsilon+\tilde{u}_{\varepsilon_{0}}|^{p-2}|\tilde{u}_{\varepsilon_{0}}||U|
+|\tilde{u}_{\varepsilon_{0}}||U|\psi_4|\Big)dx\notag\\
&\  \ \ \ \ \  \ \ \ \ \leq c\eta\Big(\|\tilde{u}_\varepsilon\|_p^p+\|\tilde{u}_{\varepsilon_{0}}\|_p^{p}
+\|U\|^p_p+\|\psi_4\|_{\frac{p}{p-2}}^{\frac{p}{p-2}}\Big).
\end{align}
By  (\ref{a2}) and (\ref{U02}),connection with $\psi_2\in L^q$, the third term on the right hand side of ({\ref{U04}}) is bounded by
\begin{align}\label{U06}
(e^{-\varepsilon\omega(t)}-e^{-\varepsilon_{0}\omega(t)})\int\limits_{\mathbb{R}^N}f(x,\tilde{u}_{\varepsilon_{0}})Udx
&\leq |e^{-\varepsilon\omega(t)}-e^{-\varepsilon_{0}\omega(t)}|\int\limits_{\mathbb{R}^N}(\varepsilon_2|\tilde{u}_{\varepsilon_{0}}|^{p-1}+\psi_2)|U|dx\notag\\
&\leq c\eta\Big(\|\tilde{u}_{\varepsilon_{0}}\|_p^{p}+\|U\|_p^{p}+\|\psi_2\|_q^q\Big),
\end{align} where $q=\frac{p}{p-1}$.
Then combination (\ref{U04})-(\ref{U06}), we find that
for every $\varepsilon\in (\varepsilon_{0}-\chi,\varepsilon_{0}+\chi)$ and $t\in [\tau,\tau+T]$,
\begin{align}\label{U07}
&\ \  \ \ \ \  \int\limits_{\mathbb{R}^N}\Big(e^{-\varepsilon\omega(t)}f(x,\tilde{u}_\varepsilon)-e^{-\varepsilon_{0}\omega(t)}f(x,\tilde{u}_{\varepsilon_{0}})\Big)Udx\notag\\
&\  \ \ \ \ \ \ \ \ \ \ \ \ \ \ \leq \alpha_3\|U\|^2+c\eta+c\eta \Big(\|\tilde{u}_\varepsilon\|_p^p+\|\tilde{u}_{\varepsilon_{0}}\|_p^{p}+\|U\|^p_p\Big)\notag\\
&\  \ \ \ \ \ \ \ \ \ \ \ \ \ \ \leq \alpha_3\|U\|^2+c_0\eta \Big(\|\tilde{u}_\varepsilon\|_p^p+\|\tilde{u}_{\varepsilon_{0}}\|_p^{p}+\|U\|^p_p\Big).
\end{align}
For the last two terms on the right hand side of ({\ref{U04}}),  by (\ref{U02}),
 we have for every $\varepsilon\in (\varepsilon_{0}-\chi,\varepsilon_{0}+\chi)$ and $t\in [\tau,\tau+T]$,
 \begin{align}\label{U08}
(e^{-\varepsilon\omega(t)}-
e^{-\varepsilon_{0}\omega(t)})\int\limits_{\mathbb{R}^N} g(t,x)Udx\leq \eta\|U\|^2+\eta\|g(t,.)\|^2,
\end{align}
\begin{align}\label{U081}
(e^{-\varepsilon\omega(t)}-
e^{-\varepsilon_{0}\omega(t)})\int\limits_{\mathbb{R}^N} h(t,x)Vdx\leq \eta\|V\|^2+\eta\|h(t,.)\|^2.
\end{align}
Then by ({\ref{U03}}) and ({\ref{U07}})-({\ref{U081}}), we get that for every $\varepsilon\in (\varepsilon_{0}-\chi,\varepsilon_{0}+\chi)$ and $t\in [\tau,\tau+T]$,
\begin{align}\label{U09}
& \ \ \ \ \ \ \ \ \ \ \ \ \ \ \ \frac{d}{dt}(\beta\|U\|^2+\alpha\|V\|^2)\leq c_1(\beta\|U\|^2+\alpha\|V\|^2)\notag\\
&+c_2\eta \Big(\|\tilde{u}_\varepsilon\|_p^p+\|\tilde{u}_{\varepsilon_{0}}\|_p^{p}
+\|u_{\varepsilon_{0}}\|_p^p+\||u_{\varepsilon}\|^p_p+\|g(t,.)\|^2+\|h(t,.)\|^2\Big),
\end{align} where $c_1$ and $c_2$ are positive constants independent of $\tau,\omega$ and $\varepsilon$.
By ({\ref{U09}}) we immediately have for every $\varepsilon\in (\varepsilon_{0}-\chi,\varepsilon_{0}+\chi)$ and $t\in [\tau,\tau+T]$,
\begin{align}\label{U10}
& \ \ \ \ \ \ \ \ \ \ \ \ \ \ \ \ \  \ \|U(t)\|^2+\|V(t)\|^2\leq c_3e^{c_1(t-\tau)}(\|U(\tau)\|^2+\|V(\tau)\|^2)\notag\\
&+c_4\eta e^{c_1(t-\tau)}\int\limits_{\tau}^t\Big(\|\tilde{u}_\varepsilon(s)\|_p^p+\|\tilde{u}_{\varepsilon_{0}}(s)\|_p^{p}
+\||u_{\varepsilon}(s)\|^p_p+\|u_{\varepsilon_{0}}(s)\|_p^p+\|g(s,.)\|^2+\|h(s,.)\|^2\Big)ds.
\end{align}
Since $e^{-\varepsilon \omega(s)}$ is continuous on $\mathbb{R}$, then for every fixed
$\tau\in \mathbb{R}$ and $\omega\in \Omega$, and $s\in[\tau,\tau+T]$, there exist $\mu=\mu(\tau,\omega,T)$ and $\nu=\nu(\tau,\omega,T)$
such that for all $\varepsilon\in (0,a]$, $\mu\leq z_\varepsilon(s,\omega)\leq\nu$ for all $s\in [\tau,\tau+T]$.
 Therefor by Lemma 4.7, it follows that for all $t\in [\tau,\tau+T]$,
\begin{align}\label{U11}
  &\int\limits_{\tau}^t \|\tilde{u}_\varepsilon(s,\tau,\omega,u_\varepsilon(\tau))\|^p_pds\notag\\
  &\leq  \mu^{-2}e^{\delta(t-\tau)}\int\limits_{\tau}^t e^{\delta(s-t)}z_\varepsilon^2(s,\omega)\|\tilde{u}_\varepsilon(s,\tau,\omega,u_\varepsilon(\tau))\|^p_pds\notag\\
 &\leq \mu^{-2}e^{\delta(t-\tau)} \Big(z_\varepsilon^2(\tau,\omega)(\|\tilde{u}_\varepsilon(\tau)\|^2+\|\tilde{v}_\varepsilon(\tau)\|^2)+c\nu^2\int\limits_{\tau}^t(\|g(s,.)\|^2+\|h(s,.)\|^2+1)ds\Big)\notag\\
 &\leq \mu^{-2}e^{\delta(t-\tau)} \Big(e^{2a|\omega(\tau)|}(\|\tilde{u}_\varepsilon(\tau)\|^2+\|\tilde{v}_\varepsilon(\tau)\|^2)+c\nu^2\int\limits_{\tau}^t(\|g(s,.)\|^2+\|h(s,.)\|^2+1)ds\Big).
\end{align}
By a similar technique we can calculate that for all $t\in [\tau,\tau+T]$,
\begin{align}\label{U12}
  &\int\limits_{\tau}^t \|u_\varepsilon(s,\tau,\omega,u_\varepsilon(\tau))\|^p_pds\notag\\
  &\leq \mu^{p-2}e^{\delta(t-\tau)} \Big(e^{2a|\omega(\tau)|}(\|\tilde{u}_\varepsilon(\tau)\|^2+\|\tilde{v}_\varepsilon(\tau)\|^2)+c\nu^2\int\limits_{\tau}^t(\|g(s,.)\|^2+\|h(s,.)\|^2+1)ds\Big).
\end{align}
Then by ({\ref{U10}}-(\ref{U12}) it gives that
\begin{align}\label{U13}
&\ \ \ \ \ \ \ \ \ \ \|U(t)\|^2+\|V(t)\|^2\leq c_3e^{c_1(t-\tau)}(\|U(\tau)\|^2+\|V(\tau)\|^2)\notag\\
&+c_5\eta e^{c_6(t-\tau)}\Big(e^{2a|\omega(\tau)|}(\mu^{-2}+\mu^{p-2})(\|\tilde{u}_\varepsilon(\tau)\|^2+\|\tilde{u}_{\varepsilon_0,\tau}\|^2+\|\tilde{v}_\varepsilon(\tau)\|^2+\|\tilde{v}_{\varepsilon_0,\tau}\|^2)\notag\\
&\ \ \ \  \ \ \ \ +\nu^2(\mu^{-2}+\mu^{p-2})\int\limits_{\tau}^t(\|g(s,.)\|^2+\|h(s,.)\|^2+1)ds+\int\limits_{\tau}^t(\|g(s,.)\|^2+\|h(s,.)\|^2)ds\Big).
\end{align}
On the other hand,  by ({\ref{U02}}) it follows that
 for every $\varepsilon\in (\varepsilon_{0}-\chi,\varepsilon_{0}+\chi)\subset(0,a]$,
 \begin{align}\label{U14}
\|U(\tau)\|^2&=\|e^{-\varepsilon\omega(\tau)}\tilde{u}_\varepsilon(\tau)-e^{-\varepsilon_{0}\omega(\tau)}\tilde{u}_{\varepsilon_0,\tau}\|^2&\notag\\
&\leq 2e^{-2\varepsilon\omega(\tau)}\|\tilde{u}_\varepsilon(\tau)-\tilde{u}_{\varepsilon_0,\tau}\|^2+
2|e^{-\varepsilon\omega(\tau)}-e^{-\varepsilon_{0}\omega(\tau)}|^2\|\tilde{u}_{\varepsilon_0,\tau}\|^2\notag\\
&\leq 2e^{2a|\omega(\tau)|}\|\tilde{u}_\varepsilon(\tau)-\tilde{u}_{\varepsilon_0,\tau}\|^2+
2\eta^2\|\tilde{u}_{\varepsilon_0,\tau}\|^2.
\end{align}
Similarly,
\begin{align}\label{U141}
\|V(\tau)\|^2\leq 2e^{2a|\omega(\tau)|}\|\tilde{v}_\varepsilon(\tau)-\tilde{v}_{\varepsilon_0,\tau}\|^2+
2\eta^2\|\tilde{v}_{\varepsilon_0,\tau}\|^2.
\end{align}
We now let $\varepsilon\rightarrow\varepsilon_{0}$ and $\|u_{\varepsilon,\tau}-u_{\varepsilon_0,\tau}\|\rightarrow0$. Then by ({\ref{U13}})-({\ref{U141}}) we obtain that for all $t\in [\tau,\tau+T]$,
\begin{align}\label{U15}
\|U(t)\|^2+\|V(t)\|^2&=\|u_{\varepsilon}(t,\tau,\omega,u_{\varepsilon,\tau})-u_{\varepsilon_{0}}(t,\tau,\omega,u_{\varepsilon_0,\tau})\|^2\notag\\
&+\|v_{\varepsilon}(t,\tau,\omega,v_{\varepsilon,\tau})-v_{\varepsilon_{0}}(t,\tau,\omega,v_{\varepsilon_0,\tau})\|^2\rightarrow0.
\end{align}
Notice that by (\ref{U02}) we also have for every $\varepsilon\in (\varepsilon_{0}-\chi,\varepsilon_{0}+\chi)$ and $t\in [\tau,\tau+T]$,
\begin{align}\label{U155}
\|\tilde{u}_\varepsilon(t,\tau,\omega ,\tilde{u}_\varepsilon(\tau))-\tilde{u}_{\varepsilon_{0}}(t,\tau,\omega ,\tilde{u}_{\varepsilon_0,\tau})\|^2
&\leq 2e^{2a|\omega(t)|}\|u_\varepsilon(t)-u_{\varepsilon_{0}}(t)\|^2+
2\eta^2\|u_{\varepsilon_{0}}(t)\|^2,
\end{align}
\begin{align}\label{U1551}
\|\tilde{v}_\varepsilon(t,\tau,\omega ,\tilde{v}_\varepsilon(\tau))-\tilde{v}_{\varepsilon_{0}}(t,\tau,\omega ,\tilde{v}_{\varepsilon_0,\tau})\|^2
&\leq 2e^{2a|\omega(t)|}\|v_\varepsilon(t)-v_{\varepsilon_{0}}(t)\|^2+
2\eta^2\|v_{\varepsilon_{0}}(t)\|^2.
\end{align}
Then by (\ref{U15})-(\ref{U1551}) we get (\ref{sss1}). Repeating the same arguments we can derive (\ref{sss2}).  $\ \ \  \ \ \ \ \ \ \Box$\\

\subsection{Main results}

We are now at the point to present the main results in this paper.\\

\textbf{Theorem 4.9.} \emph{Suppose $\varepsilon\in\mathbb{R}$ and (\ref{a1})-(\ref{a5}) hold true. Then}

\emph{(i)\ \ Each random cocycle $\varphi_\varepsilon$ generated by
(\ref{FN1})-(\ref{FN3}) has a unique pullback attractor $\mathcal{A}_\varepsilon$ and the corresponding deterministic cocycle $\varphi_0$ has a unique pullback attractor
$\mathcal{A}_0$ in $L^2(\mathbb{R}^N)\times L^2(\mathbb{R}^N)$, both $\mathcal{A}_\varepsilon$ and $\mathcal{A}_0$ are $(L^2(\mathbb{R}^N)\times L^2(\mathbb{R}^N),L^l(\mathbb{R}^N)\times L^2(\mathbb{R}^N))$-pullback attractors.}

\emph{(ii) If further (\ref{a6}) holds true, then the family $\mathcal{A}_\varepsilon$ is upper semi-continuous under the
Hausdorff semi-distance of $L^l(\mathbb{R}^N)\times L^2(\mathbb{R}^N)$ at any $\varepsilon_0\in\mathbb{R}$. Here $l\in (2,p],p>2$.}
\\

\emph{Proof}\ \   Let $X=L^2(\mathbb{R}^N)\times L^2(\mathbb{R}^N)$ and $L^l(\mathbb{R}^N)\times L^2(\mathbb{R}^N)$. Then it is known that
the hypothesises (\emph{H1})) and (\emph{H2})) (see Lemma 2.7 in \cite{Zhao3}) hold true.
By the Sobolev interpolation and association with Lemma 4.2  and Lemma 4.6, we immediately
 obtain the uniformly pullback asymptotic compactness in $L^l(\mathbb{R}^N)\times L^2(\mathbb{R}^N)$ for any $2<l<p$.
Then along with  uniform absorption (see Lemma 4.1) and convergence property (see Lemma 4.8), all conditions of Theorem 2.10 are satisfied.
$\ \ \  \ \ \ \ \ \ \Box$\\

\section{Existence of random equilibria for the generated random cocycle}

Random equilibrium is  a special case of omega-limit sets.
We can refer to \cite{Arn,Chues} for the definitions and applications to monotone random dynamical system.
The problem of the construction of equilibria for a general random dynamical system is rather complicate \cite{Chues}. Recently, Gu \cite{Gu1} proved that the stochastic FitzHugh-Nagumo lattice equations driven by fractional Brownian motions possess a unique equilibrium.
However, we here introduce the random equilibrium  in the case of non-autonomous stochastic dynamical system. Specifically, we have\\

\textbf{Definition 5.1.} \emph{Let $(Q,\{\sigma_t\}_{\mathbb{R}})$ and $(\Omega,\mathcal {F},{P},\{\vartheta_t\}_{t\in\mathbb{R}})$ be parametric dynamical systems. A random variable $u^*: Q\times\Omega\mapsto X$
is said to be an equilibrium (or fixed point, or stationary solution) of random  cocycle $\varphi$ if it is invariant under $\varphi$, i.e., if
$$
\varphi(t,q,\omega, u^*(q,\omega))=u^*(\sigma_tq,\vartheta_t\omega)\ \  for\ all\ t\geq0,\ q \in Q,\ \omega\in\Omega.
$$}

In this section, the parametric dynamical systems $(Q,\{\sigma_t\}_{\mathbb{R}})$ and $(\Omega,\mathcal {F},{P},\{\vartheta_t\}_{t\in\mathbb{R}})$
are same as in section 3. We will prove the existence of equilibrium
for problem (\ref{FN1})-(\ref{FN3}) on the whole space $\mathbb{R}^N$. To this end, we need to assume that
\begin{align} \label{AAA}
 \delta=\min\{\lambda,\sigma\}>\alpha_3,\ \ \  \beta\geq1,
\end{align} where $\alpha_3$ is as in (\ref{a3}) and $\lambda,\sigma,\beta$ are as in the FitzHugh-Nagumo system (\ref{FN1})-(\ref{FN2}).

 For convenience, here we write $\varepsilon=1$.
 First, we have\\

\textbf{Lemma 5.2.} \emph{Suppose that $g$ and $h$ satisfy (\ref{a5}) and $f$ satisfies (\ref{a1}) and (\ref{a3}) such that (\ref{AAA}) holds. Then there exists a positive constant $0<b_0<\delta-\alpha_3$ such that
the solutions of problem (\ref{FN1})-(\ref{FN3}) with initial values $(\tilde{u}_{\tau-t_i}, \tilde{v}_{\tau-t_i})(i=1,2), t_1<t_2$ satisfy the following decay  property:
\begin{align} \label{}
&\|\tilde{u}(\tau, \tau-t_1,\vartheta_{-\tau}\omega, \tilde{u}_{\tau-t_1})-\tilde{u}(\tau, \tau-t_2,\vartheta_{-\tau}\omega, \tilde{u}_{\tau-t_2})\|^2\notag\\
&+\|\tilde{v}(\tau, \tau-t_1,\vartheta_{-\tau}\omega, \tilde{v}_{\tau-t_1})-\tilde{v}(\tau, \tau-t_2,\vartheta_{-\tau}\omega, \tilde{v}_{\tau-t_2})\|^2\notag\\
&\leq c\Big(e^{-b_0t_1}z^2(-t_1,\omega)(\|\tilde{u}_{\tau-t_1}\|^2+\|\tilde{v}_{\tau-t_1}\|^2)+e^{-b_0 t_2}z^2(-t_2,\omega)(\|\tilde{u}_{\tau-t_2}\|^2+\|\tilde{v}_{\tau-t_2}\|^2)\Big)\notag\\
 &+ce^{(b_0-b)t_1}\int_{-\infty}^{0}e^{b_0 s}z^2(s,\omega)(\|g(s+\tau,.)\|^2+\|h(s+\tau,.)\|^2+1)ds,\notag
\end{align}where $c$ is a deterministic non-random constant.}\\

\emph{Proof}\ \  Put $\bar{u}=u(t, \tau-t_1,\vartheta_{-\tau}\omega, u_{\tau-t_1})-u(t, \tau-t_2,\vartheta_{-\tau}\omega, u_{\tau-t_2})$ and
 $$\bar{v}=v(t, \tau-t_1,\vartheta_{-\tau}\omega, v_{\tau-t_1})-v(t, \tau-t_2,\vartheta_{-\tau}\omega, v_{\tau-t_2}).$$
  Then from (\ref{pr1})-(\ref{pr2}), along with (5.1),  we have
\begin{align} \label{equ1}
 \frac{d}{dt}(\beta\|\bar{u}\|^2+\alpha\|\bar{v}\|^2)+b(\beta\|\bar{u}\|^2+\alpha\|\bar{v}\|^2)\leq 0,
\end{align} where $b=\delta-\alpha_3$.
By applying Gronwall lemma to (\ref{equ1}) over the interval $[\tau-t_1,\tau]$, we immediately get
\begin{align} \label{equ2}
\|\bar{u}(\tau)\|^2+\|\bar{v}(\tau)\|^2 &\leq ce^{-bt_1}(\|u(\tau-t_1, \tau-t_2,\vartheta_{-\tau}\omega, u_{\tau-t_2})-u_{\tau-t_1}\|^2\notag\\
&+\|v(\tau-t_1, \tau-t_2,\vartheta_{-\tau}\omega, v_{\tau-t_2})-v_{\tau-t_1}\|^2)\notag\\
 &\leq ce^{-bt_1}(\|u(\tau-t_1, \tau-t_2,\vartheta_{-\tau}\omega, u_{\tau-t_2})\|^2\notag\\
 &+\|v(\tau-t_1, \tau-t_2,\vartheta_{-\tau}\omega, v_{\tau-t_2})\|^2)+ce^{-bt_1}(\|u_{\tau-t_1}\|^2+\|v_{\tau-t_1}\|^2),
\end{align} where $c=c(\alpha,\beta)$ is a positive  deterministic constant.
Write
\begin{align} \label{equ3}
0<\delta_0<b_0<b=\delta-\alpha_3,
\end{align} where $\delta_0$ is as in (\ref{a5}).
From (\ref{energy}) and using (\ref{equ3}),  we have
\begin{align} \label{equ4}
\frac{d}{dt}(\beta\|u\|^2+\alpha\|v\|^2)&+b_0(\beta\|u\|^2+\alpha\|v\|^2)\leq cz^2(t,\omega)(\|g(t,.)\|^2+\|h(t,.)\|^2+1).
\end{align}
Then by Gronwall lemma again over the interval $[\tau-t_2,\tau-t_1]$, we find that
\begin{align} \label{equ44}
 &\|u(\tau-t_1, \tau-t_2,\vartheta_{-\tau}\omega, u_{\tau-t_2})\|^2+\|v(\tau-t_1, \tau-t_2,\vartheta_{-\tau}\omega, v_{\tau-t_2})\|^2\notag\\
 &\ \ \ \ \ \leq ce^{b_0(t_1-t_2)}(\|u_{\tau-t_2}\|^2+\|v_{\tau-t_2}\|^2)\notag\\
 &\ \ \ \  \ \ \  +c\int\limits_{\tau-t_2}^{\tau-t_1}e^{-b_0(\tau-t_1-s)}z^2(s,\vartheta_{-\tau}\omega)(\|g(s,.)\|^2+\|h(s,.)\|^2+1)ds\notag\\
 & \ \  \ \ \ \ \ \ \ \leq  ce^{b_0(t_1-t_2)}(\|u_{\tau-t_2}\|^2+\|v_{\tau-t_2}\|^2)\notag\\
 & \ \ \  \ \ \ \ \ \ \ \ \ +ce^{b_0 t_1}e^{2\omega(-\tau)}\int\limits_{-\infty}^{0}e^{b_0 s}z^2(s,\omega)(\|g(s+\tau,.)\|^2+\|h(s+\tau,.)\|^2+1)ds,
\end{align} where $c=c(\alpha,\beta)$ is a positive  deterministic constant. Then combination  (\ref{equ44}) and (\ref{equ2}) we have
\begin{align} \label{equ5}
 \|\bar{u}(\tau)\|^2+\|\bar{v}(\tau)\|^2
 &\leq ce^{-bt_1}(\|u_{\tau-t_1}\|^2+\|v_{\tau-t_1}\|^2)+ce^{(b_0-b)t_1}e^{-b_0 t_2}(\|u_{\tau-t_2}\|^2+\|v_{\tau-t_2}\|^2)\notag\\
 &+ce^{(b_0-b)t_1}e^{2\omega(-\tau)}\int_{-\infty}^{0}e^{b_0 s}z^2(s,\omega)(\|g(s+\tau,.)\|^2+\|h(s+\tau,.)\|^2+1)ds\notag\\
 &\leq ce^{-b_0t_1}(\|u_{\tau-t_1}\|^2+\|v_{\tau-t_1}\|^2)+e^{-b_0 t_2}(\|u_{\tau-t_2}\|^2+\|v_{\tau-t_2}\|^2)\notag\\
 &+ce^{(b_0-b)t_1}e^{2\omega(-\tau)}\int_{-\infty}^{0}e^{b_0 s}z^2(s,\omega)(\|g(s+\tau,.)\|^2+\|h(s+\tau,.)\|^2+1)ds,
\end{align}
where we have used $e^{(b_0-b)t_1}\leq 1$ for $b_0<b$. In terms of the relation (\ref{trans}),  we get
\begin{align} \label{}
 \|\bar{\tilde{u}}(\tau)\|^2+\|\bar{\tilde{v}}(\tau)\|^2&\leq
 e^{-2\omega(-\tau)}\|\bar{u}(\tau)\|^2+\|\bar{v}(\tau)\|^2\notag\\
 &\leq ce^{-2\omega(-\tau)} \Big(e^{-b_0t_1}(\|u_{\tau-t_1}\|^2+\|v_{\tau-t_1}\|^2)+e^{-b_0 t_2}(\|u_{\tau-t_2}\|^2+\|v_{\tau-t_2}\|^2)\Big)\notag\\
 &+ce^{(b_0-b)t_1}\int_{-\infty}^{0}e^{b_0 s}z^2(s,\omega)(\|g(s+\tau,.)\|^2+\|h(s+\tau,.)\|^2+1)ds\notag\\
 &=ce^{-2\omega(-\tau)} \Big(e^{-b_0t_1}z^2(\tau-t_1,\vartheta_{-\tau}\omega)(\|\tilde{u}_{\tau-t_1}\|^2+\|\tilde{v}_{\tau-t_1}\|^2)\notag\\
 &+e^{-b_0 t_2}z^2(\tau-t_2,\vartheta_{-\tau}\omega)(\|\tilde{u}_{\tau-t_2}\|^2+\|\tilde{v}_{\tau-t_2}\|^2)\Big)\notag\\
 &+ce^{(b_0-b)t_1}\int_{-\infty}^{0}e^{b_0 s}z^2(s,\omega)(\|g(s+\tau,.)\|^2+\|h(s+\tau,.)\|^2+1)ds\notag\\
 &=c\Big(e^{-b_0t_1}z^2(-t_1,\omega)(\|\tilde{u}_{\tau-t_1}\|^2+\|\tilde{v}_{\tau-t_1}\|^2)+e^{-b_0 t_2}z^2(-t_2,\omega)(\|\tilde{u}_{\tau-t_2}\|^2+\|\tilde{v}_{\tau-t_2}\|^2)\Big)\notag\\
 &+ce^{(b_0-b)t_1}\int_{-\infty}^{0}e^{b_0 s}z^2(s,\omega)(\|g(s+\tau,.)\|^2+\|h(s+\tau,.)\|^2+1)ds,\notag
\end{align} which finishes the proof.$\ \ \  \ \ \ \ \ \ \Box$\\

\textbf{Lemma 5.3.} \emph{Suppose that $g$ and $h$ satisfy (\ref{a5}), $f$ satisfies (\ref{a1}) and (\ref{a3}) such that (\ref{AAA}) holds. Let $ D=\{D(\tau,\omega);\tau\in\mathbb{R},\omega\in\Omega\}\in\mathcal{D}_\delta$ where $\mathcal{D}_\delta$ is defined as in (\ref{D}). Then for $\tau\in\mathbb{R},\omega\in\Omega$, there exists  a unique element $u^*=u^*(\tau,\omega)\in L^2(\mathbb{R}^N)\times L^2(\mathbb{R}^N)$ such that
\begin{align} \label{}
\lim\limits_{t\rightarrow+\infty}(\tilde{u}(\tau, \tau-t,\vartheta_{-\tau}\omega, \tilde{u}_{\tau-t}),\tilde{u}(\tau, \tau-t,\vartheta_{-\tau}\omega, \tilde{u}_{\tau-t}))=u^*(\tau,\omega),\ \  \mbox{in}\ L^2(\mathbb{R}^N)\times L^2(\mathbb{R}^N),\notag
\end{align} where $(\tilde{u}_{\tau-t},\tilde{v}_{\tau-t})\in D(\tau-t,\vartheta_{-t}\omega)$.
Furthermore, the convergence is uniform (w.r.t. $(\tilde{u}_{\tau-t},\tilde{v}_{\tau-t})\in D(\tau-t,\vartheta_{-t}\omega)$).}\\

\emph{Proof}\ \ We choose the constant $b_0$ in Lemma 5.2 satisfying $b>b_0>\delta_1$, where $\delta_1$ is as in (\ref{D}).
 If $(\tilde{u}_{\tau-t_i},\tilde{v}_{\tau-t_i})\in D(\tau-t_i,\vartheta_{-t_i}\omega)$, $i=1,2$, then similar to (\ref{un02}) we have
 $$\lim\limits_{t_i\rightarrow+\infty}e^{-b_0t_i}z^2(-t_i,\omega)\|(\tilde{u}_{\tau-t_i},\tilde{v}_{\tau-t_i})\|^2=0,i=1,2.$$
  Thus the result is derived directly from Lemma 5.2.
 $\ \ \  \ \ \ \ \ \ \Box$\\

\textbf{Lemma 5.4.} \emph{Suppose that $g$ and $h$ satisfies (\ref{a5}), $f$ satisfies (\ref{a1}) and (\ref{a3}) such that (\ref{AAA}) holds. Then for $\tau\in\mathbb{R},\omega\in\Omega$, the  element $u^*=u^*(\tau,\omega)$ defined in Lemma 5.3 is a unique random equilibrium for the cocycle $\varphi$  defined by (\ref{eq0}) in $L^2(\mathbb{R}^N)\times L^2(\mathbb{R}^N) $, i.e.,
\begin{align} \label{}
\varphi(t, \tau,\omega, u^*(\tau,\omega))=u^*(\tau+t,\vartheta_{t}\omega),\ \ \mbox{for every }\ t\geq0,\ \tau\in\mathbb{R},\ \omega\in\Omega.\notag
\end{align}
Furthermore, the random equilibrium $\{u^*(\tau,\omega), \tau\in\mathbb{R},\ \omega\in\Omega\}$  is the unique element of the pullback attractor $\mathcal{A}=\{\mathcal{A}(\tau,\omega); \tau\in\mathbb{R},\ \omega\in\Omega\}$ for the random cocycle $\varphi$, i.e., for every $\tau\in\mathbb{R},\omega\in\Omega$,
$\mathcal{A}(\tau,\omega)=\{u^*(\tau,\omega)\}.$}\\

\emph{Proof}\ \  From the definition of random cocycle,
$$
\varphi(t,\tau-t, \vartheta_{-t}\omega, (\tilde{u}_{\tau-t},\tilde{v}_{\tau-t}))=(\tilde{u}(\tau,\tau-t,\vartheta_{-\tau}\omega,\tilde{u}_{\tau-t}),\tilde{v\emph{}}(\tau,\tau-t,\vartheta_{-\tau}\omega,\tilde{v}_{\tau-t}) ),
$$
then for for every $\tau\in\mathbb{R},\omega\in\Omega$, we see from Lemma 5.3 that
\begin{align} \label{6.8}
u^*(\tau,\omega)=\lim\limits_{t\rightarrow+\infty} \varphi(t, \tau-t,\vartheta_{-t}\omega, (\tilde{u}_{\tau-t},\tilde{v}_{\tau-t})),
\end{align} where $(\tilde{u}_{\tau-t},\tilde{v}_{\tau-t})\in D(\tau-t,\vartheta_{-t}\omega)$. Thus by the continuity and the cocycle property  of $\varphi$ and (\ref{6.8}), we find that for every $t\geq0,\tau\in\mathbb{R},\omega\in\Omega$,
\begin{align} \label{}
\varphi(t, \tau,\omega, u^*(\tau,\omega))&=\varphi(t, \tau,\omega, .)\circ \lim\limits_{s\rightarrow+\infty} \varphi(s, \tau-s,\vartheta_{-s}\omega, (\tilde{u}_{\tau-t},\tilde{v}_{\tau-t}))\notag\\
&=\lim\limits_{s\rightarrow+\infty} \varphi(t, \tau,\omega, .)\circ \varphi(s, \tau-s,\vartheta_{-s}\omega, (\tilde{u}_{\tau-t},\tilde{v}_{\tau-t}))\notag\\
&=\lim\limits_{s\rightarrow+\infty} \varphi(t+s, \tau-s,\vartheta_{-s}\omega, (\tilde{u}_{\tau-t},\tilde{v}_{\tau-t}))\notag\\
&=\lim\limits_{s\rightarrow+\infty} \varphi(t+s, (\tau+t)-t-s, \vartheta_{-s-t}\vartheta_{t}\omega, (\tilde{u}_{\tau-t},\tilde{v}_{\tau-t}))=u^*(\tau+t,\vartheta_{t}\omega),\notag
\end{align}
which also implies  the invariance of $\mathcal{A}$, that is, $\varphi(t, \tau,\omega, \mathcal{A}(\tau,\omega))=\mathcal{A}(\tau+t,\vartheta_t\omega)$. The compactness of
$\mathcal{A}(\tau,\omega)$ is obvious and the attracting property follows from (\ref{6.8}). $\ \ \  \ \ \ \ \ \ \Box$\\

\textbf{Acknowledgments:}

 This work was supported by  Chongqing Basis and Frontier Research Project of China (No. cstc2014jcyjA00035), National Natural Science Foundation of China
(No. 11271388), China Postdoctoral Science Foundation Grant (No.2014M550452).
 \\

\end{document}